\documentclass[11pt,a4paper]{article}

\usepackage{graphicx,latexsym}
\usepackage{color}
\usepackage{keyval}
\usepackage{amssymb,amsmath}
\usepackage{amsfonts}
\usepackage{amsthm}
\usepackage{dsfont}
\usepackage{mathrsfs}
\usepackage{cite}
\usepackage{geometry}

\newtheorem{remark}{\textbf{Remark}}[section]

\allowdisplaybreaks[4]
\def\d{\mathrm{d}}

\begin{document}
\title{
A Review on the Cahn--Hilliard Equation: Classical Results and Recent Advances in Dynamic Boundary Conditions
}
\author{
{\sc Hao Wu} \\ School of Mathematical Sciences\\
Fudan University\\
Handan Road 220, Shanghai 200433, CHINA\\
Email: \texttt{haowufd@fudan.edu.cn}
}
\date{}
\maketitle

\begin{abstract}
The Cahn--Hilliard equation is a fundamental model that describes the phase separation process in multi-component mixtures. It has been successfully extended to many different contexts in several scientific fields. In this survey article, we briefly review the derivation, structure as well as some analytical issues for the Cahn--Hilliard equation and its variants. Our focus will be placed on the well-posedness and long-time behavior of the Cahn--Hilliard equation in the classical setting and recent progresses on the dynamic boundary conditions accounting for non-trivial boundary effects.


\end{abstract}

\section{Introduction}

The Cahn--Hilliard equation provides a continuous description of the phase separation process for binary mixtures. It was first proposed in \cite{CH58,CH61} to model the so-called spinodal decomposition of binary alloys in a rapid cooling process, assuming isotropy of the material. As pointed out in \cite{NC08}, the Cahn--Hillirad equation is rather broad ranged in its evolution scope such that it is capable of describing important qualitative features of many systems undergoing phase separation at different time stages. Besides the spontaneous heterogenization of a binary mixture like spinodal decomposition, it can also model other mechanisms in pattern formation, for instance, the process of nucleation and growth, and the process of coarsening \cite{BF93,CH59,NC08}. The Cahn--Hilliard equation is a representative of the so-called diffuse interface models describing the evolution of free interfaces during phase transitions. Instead of the classical sharp-interface formulation, the diffuse interface models represent free interfaces that separate different components of the mixture as thin layers with finite thickness over which material properties vary smoothly \cite{AN97}. This approach has several advantages. For example, explicit tracking of the interface, which is usually a difficult task, can be avoided in both mathematical formulation and numerical computations. Besides, evolution of complex geometries and topological changes of the interface can be handled in a natural way, see \cite{DF20} for a recent review. In the past years, the Cahn--Hilliard equation and its variants have been successfully applied in different fields and turned out to be efficient tools for the study of a wide variety of segregation-like phenomena \cite{Kim16}, for instance, diblock copolymer \cite{OK86}, image inpainting \cite{BEG07,BEG07a}, tumor growth \cite{GLSS16,Od10,OP13}, biology \cite{DLRW,KS08}, two-phase flows \cite{AGG12,AN97,Gur96a,HH77,LaW18,LS03}, moving contact line \cite{Ja00,Q1}, and so on.

Assume that $T\in (0,+\infty]$, $\Omega\subset \mathbb{R}^d$ ($d=1,2,3$) is a bounded domain with a sufficiently smooth boundary $\partial \Omega$. When $d=1$, $\Omega$ is simply an interval, e.g., $\Omega=(0,l)$ for some $l>0$. The classical Cahn--Hilliard equation can be written in the following form:
\begin{equation}
\left\{
\begin{aligned}
&\partial_t \phi=\nabla \cdot[M(\phi)\nabla \mu],\quad &\text{in}\ \Omega \times (0,T),\\
&\mu=-\epsilon^2\Delta \phi + F'(\phi),\quad &\text{in}\ \Omega \times (0,T).
\end{aligned}
\right.
\label{claCH}
\end{equation}
The unknown function $\phi$ is called the order parameter or phase-field, which is related to local concentrations of the two components of a binary mixture. It may have different interpretations according to physical contexts, for instance, the volume fraction, mass fraction, or mole fraction \cite{DF20,NC08}. Usually, people consider a rescaled form of $\phi$ such that it denotes the difference between concentrations for the two components, that is, $\phi=c_A-c_B$. Since the concentrations satisfy $c_A, c_B\in[0,1]$ and $c_A+c_B=1$, it is straightforward to check that $\phi$ should take its value in the physical interval $[-1,1]$, with $\pm 1$ corresponding to the pure states. In a general framework for the description of free interfaces, the phase-field function $\phi$ takes distinct values in different bulk phases away from the diffuse interface separating them, and the interface can be identified with an intermediate level set (e.g., the zero-level set) of $\phi$.

In equation \eqref{claCH}, $\mu$ denotes the chemical potential, $\epsilon$ is a positive constant and $M$ is a nonnegative quantity standing for the mobility that can be chosen as either a positive constant or a concentration dependent function $M=M(\phi)$. From the mathematical point of view, \eqref{claCH} is a fourth order parabolic equation for the unknown variable $\phi$ when $M>0$.  Thus, suitable initial and boundary conditions should be taken into account to form a well-posed problem. For the initial condition, we take
\begin{equation}
\phi|_{t=0}=\phi_0(x),\quad \text{in}\ \Omega. \label{ini0}
\end{equation}
On the other hand, one of the classical choices for boundary conditions is the following homogeneous Neumann type:
 \begin{align}
&M(\phi)\nabla \mu \cdot \mathbf{n}=0, \quad\  \text{on}\ \partial\Omega \times (0,T),\label{bdN1}\\
&\partial_\mathbf{n}\phi=0, \qquad\qquad \quad \text{on}\ \partial \Omega \times (0,T).\label{bdN2}
\end{align}
Here, $\mathbf{n}=\mathbf{n}(x)$ denotes the unit exterior normal to the boundary $\partial\Omega$ and $\partial_\mathbf{n}$ stands for the outward normal derivative on $\partial \Omega$. Other boundary conditions that are of interest are Dirichlet boundary conditions (for $\phi$ and $\mu$) and periodic boundary conditions (cf. \cite{E89}).

The Cahn--Hilliard equation has been extensively studied in the literature both analytically and numerically. Our aim in this paper is to review the derivation, mathematical structure and some analytical issues for the Cahn--Hilliard equation and its variants, with main focus on the well-posedness and long-time behavior of global solutions. Some other important and interesting problems like the nonlocal interactions, limit motions, the formal and rigorous justifications of the sharp-interface limit etc are beyond the scope of this review.
In Section 2, we present results for the Cahn--Hilliard equation in the classical setting, i.e., the initial boundary value problem \eqref{claCH}--\eqref{bdN2}. In Section 3, we discuss some recent developments on the dynamic boundary conditions.

\section{Cahn--Hilliard Equation in the Classical Setting}
\setcounter{equation}{0}

\subsection{Derivation and structure}

For an isotropic binary mixture with nonuniform composition at a fixed temperature,
we consider the following Ginzburg--Landau type free energy:
\begin{align}
E^\mathrm{bulk}(\phi)=\int_\Omega \frac{\epsilon^2}{2} |\nabla \phi|^2 +F(\phi) \, \d x.\label{Ebulk}
\end{align}
 The energy functional $E^\mathrm{bulk}$ consists of two contributions, that is, the bulk part and the gradient part. The bulk energy part represents the interaction of different components in a homogeneous system, where $F(\phi)$ denotes the (Helmholtz) free energy density of mixing. A typical thermodynamically relevant example is the following logarithmic potential (see \cite{CH58,Kim14})
\begin{align}
&F_{\mathrm{log}}(s)=\frac{\theta}{2}\Big[(1+s)\ln (1+s)+(1-s)\ln(1-s)\Big]-\frac{\theta_c}{2}s^2,\quad s\in(-1,1),
\label{logpo}
\end{align}
for some constants $\theta,\theta_c>0$, where $\theta$ is the absolute temperature of the system and $\theta_c$ is the critical temperature of phase separation. The potential \eqref{logpo} is also related to the Flory--Huggins free energy density in assessing the mutual miscibility of polymer solutions. When $0<\theta<\theta_c$, it is easy to verify that $F_{\mathrm{log}}$ has a double-well structure with two minima $\pm \phi_*\in (-1,1)$, where $\phi_*$ is the positive root of the equation $F_{\mathrm{log}}'(s)=0$ (see \cite{BE91}). This case is of physical importance since the phase separation may occur, see for instance, the phase diagram in \cite{NS84,NC08}. The interval $(-\phi^*,\phi^*)$ with $\phi^*=(1-\theta/\theta_c)^{1/2}>0$ is called the spinodal region, where $F_{\mathrm{log}}''(s)<0$. When a homogeneous state is located in $(-\phi^*,\phi^*)$, any  small fluctuation in composition will lower the free energy and yield spontaneous phase separation of the mixture towards the equilibrium compositions of two phases corresponding to $\pm \phi_*$. In the literature, when the system temperature $\theta$ is closed to the critical temperature $\theta_c$, i.e., the ``shallow quenching" case, the singular potential $F_{\mathrm{log}}$ is often approximated by a polynomial of degree four like the following form
\begin{equation}
F_{\mathrm{reg}}(s)=\frac{1}{4}\big(1-s^2\big)^2,\quad s\in\mathbb{R}. \label{regpo}
\end{equation}
More precisely, in view of \eqref{logpo}, we can apply Taylor's expansion at $s=0$ to get (cf. \cite{DD95})
$$
F_{\mathrm{log}}(s)\approx \left(\frac{\theta}{2}-\frac{\theta_c}{2}\right)s^2+\frac{\theta}{12}s^4.
$$
For the special case $\theta=3$, $\theta_c=4$, we recover \eqref{regpo} (after adding the resultant with a constant $1/4$). The above simple approximation brings great convenience in the mathematical analysis and numeric simulation for the Cahn--Hilliard equation. Nevertheless, we note that it yields a shift of the location of minima, i.e., from $\pm \phi_*$ to $\pm1$. In the literature, there is another type of singular potential called the double obstacle potential, which can be obtained in the ``deep-quench limit" of the logarithmic potential \eqref{logpo} by letting $\theta\to 0^+$ (see \cite{BE91}). In this case, the spinodal region eventually expands to $(-1,1)$ and the corresponding free energy density usually takes the following form (after a shift by the constant $\theta_c/2$):
\begin{equation}
F_{\mathrm{obs}}(s)
=\frac{\theta_c}{2}(1-s^{2})+\mathbb{I}_{[-1,1]}(s)=
\begin{cases}
\displaystyle{\frac{\theta_{c}}{2}(1-s^2),}\quad\ \text{if}\ s\in [-1,1],\\
+\infty,\qquad \qquad \text{if}\ |s|>1.
\end{cases}
\label{obspo}
\end{equation}

The gradient part of $E^\mathrm{bulk}$ represents the spatial variation in composition of the mixture, where $\epsilon^2>0$ is a parameter often called the coefficient of gradient energy \cite{CH58,NS84,E89}. The small positive constant $\epsilon$ measures the capillary width of the mixture (i.e., width of the transition layer), see \cite{AN97,LS03,DF20}. Formally speaking, the gradient term can be obtained by making expansion of a free energy density in a region of nonuniform composition $\widetilde{F}(\phi)(\phi,\nabla \phi, \nabla^2\phi,...)\approx F(\phi)+ (\epsilon^2/2) |\nabla\phi|^2+ \cdots $, see e.g.,  \cite{CH58,NS84} for details. As pointed out in \cite{Fife10}, the gradient part serves as a penalty for the phase-field function $\phi$ having sudden changes with respect to the spatial variable $x$ and corresponds to the tendency of the mixture to prefer to be uniform in space. It accounts for heterogeneity of the mixture and gives an approximation of the interfacial surface energy \cite{DF20,LS03}. We note that regularization by gradient terms can date back to the earlier work \cite{vdV1893}. Here, this gradient term can be regarded as an elliptic regularization against the possible backwards diffusion (noting that $F$ can be nonconvex in the spinodal region), in order to guarantee the PDE system to be well-posed (cf. \cite{Fife10,NC08}). \medskip

\textbf{Derivation from the mass balance law}. There are several ways to derive (and to understand) the Cahn--Hilliard equation \eqref{claCH}. We first present the derivation by using basic thermodynamics (cf. \cite{CH58,Kim14,NS84}). Recalling the definition of the order parameter $\phi$, we see that the concentrations of the two components of a binary mixture can be determined once the scalar function $\phi$ is determined. Thus, the evolution of the composition can be predicted by a single evolution equation for $\phi$.

The conserved dynamics of the phase separation process is due to the generalized (non-Fickian) diffusion driven by gradients in the chemical potential $\mu$ (see \cite{CH58}, and also \cite{E89,Kim14,Mi17,NS84}).
Phenomenologically, we consider a continuity equation
\begin{equation}
\partial_t\phi + \nabla \cdot \mathbf{J}=0,
\label{conti1}
\end{equation}
where the vector $\mathbf{J}$ denotes the mass flux.
Equation \eqref{conti1} yields a differential (local) form for the law of mass balance.
The natural and possibly the simplest boundary conditions are as follows (cf. \eqref{bdN1}--\eqref{bdN2})
 \begin{align}
&\mathbf{J} \cdot \mathbf{n}=0, \quad\  \text{on}\ \partial\Omega \times (0,T),\label{bdN1z}\\
&\partial_\mathbf{n}\phi=0, \quad \ \ \text{on}\ \partial \Omega \times (0,T).\label{bdN2z}
\end{align}
The first boundary condition \eqref{bdN1z} is usually called the no-flux boundary condition such that after integration by part, we can easily deduce from \eqref{conti1} the mass conservation of the system (at least in a formal manner):
\begin{equation*}
\frac{\d}{\d t}\int_\Omega \phi\, \d x=\int_\Omega \partial_t \phi\, \d x=-\int_\Omega \nabla \cdot \mathbf{J} \,\d x=\int_{\partial \Omega} \mathbf{J}\cdot\mathbf{n}\, \d S=0, \quad \forall\, t\in (0,T).
\label{masscov1}
\end{equation*}
Therefore, it holds
\begin{equation}
\int_\Omega \phi(x,t)\, \d x = \int_\Omega \phi_0(x)\, \d x,\quad \forall\, t\in [0,T]. \label{massz}
\end{equation}
On the other hand, the homogeneous Neumann boundary condition \eqref{bdN2z} for $\phi$ also has its physical interpretation: the free interface between the two components intersects the solid wall (i.e., the boundary $\partial \Omega$) at a perfect static contact angle of $\pi/2$ (cf. \cite{Ja00,Q1}). Under this choice, the chemical potential $\mu$, that is defined as the variational derivative of the bulk free energy $E^\mathrm{bulk}$ with respect to $\phi$, is given by
\begin{equation}
\mu=\frac{\delta E^\mathrm{bulk}(\phi)}{\delta\phi}=-\epsilon^2\Delta \phi + F'(\phi),
\label{che-mu}
\end{equation}
where $F'$ denotes the derivative of the bulk potential $F$ with respect to $\phi$ (cf. \cite{Fife10,NS84}).
The chemical potential may be viewed as a forcing term proportional to the local distance from equilibrium \cite{NS84}. The right expression of the mass flux $\mathbf{J}$ should be determined to fulfill the basic thermodynamics, that is, the evolution of $\phi$ must occur in such a way that the free energy $E^\mathrm{bulk}(\phi)$ does not increase in time. For instance, we may postulate the following constitutive equation
\begin{equation}
\mathbf{J}=-M(\phi)\nabla \mu.
\label{fluxJ}
\end{equation}
Then combining \eqref{conti1}, \eqref{che-mu} and \eqref{fluxJ}, we arrive at the Cahn--Hilliard equation \eqref{claCH}.
Besides, it is obvious that \eqref{bdN1z} together with \eqref{fluxJ} yields the boundary condition \eqref{bdN1}.
An important consequence of \eqref{bdN1z}, \eqref{bdN2z} and \eqref{fluxJ} is that the free energy $E^\mathrm{bulk}(\phi)$ is indeed non-increasing in time, which can be seen from a direct calculation (cf. \cite{E89,Kim14}):
\begin{align}
\frac{\d}{\d t} E^\mathrm{bulk}(\phi(t))
& =\int_\Omega \left(\epsilon^2 \nabla \phi \cdot \nabla \partial_t \phi
+ F'(\phi) \partial_t \phi \right)\,\d x
\notag
\\
&=\int_\Omega \mu \partial_t \phi\, \d x
 =\int_\Omega \mu \nabla \cdot (M(\phi)\nabla \mu)\,\d x\notag
\\
&=-\int_\Omega \nabla \mu\cdot(M(\phi)\nabla \mu)\,\d x +\int_{\partial\Omega} \mu(M(\phi)\nabla \mu \cdot\mathbf{n})\,\d S\notag
\\
&= -\int_\Omega M(\phi)|\nabla \mu|^2 \, \d x,\quad \forall\, t\in (0,T).
\label{BELz}
\end{align}
In this sense, \eqref{bdN2z} is sometimes referred to as the variational boundary condition such that it guarantees the validity of the energy dissipation law \eqref{BELz} and thus the resulting system fulfills the requirement from thermodynamics.

When the Cahn--Hilliard equation \eqref{claCH} is subject to other types of boundary conditions for $\phi$ and $\mu$, for instance, the periodic boundary conditions or the (nonhomogeneous) Dirichlet boundary conditions, a similar energy dissipation law can still be derived, see e.g., \cite{E89}.\medskip

\textbf{Gradient flow structure}. Another possible derivation of \eqref{claCH} we would like to mention is based on the gradient flow approach, see e.g., \cite{Fife10}, where the author treated the case with $M$ being a positive constant. Within this framework, the Cahn--Hilliard equation \eqref{claCH} can be derived by considering the constrained gradient dynamics for the free energy functional $E^\mathrm{bulk}$, subject to the mass conservation \eqref{massz}. Namely, we seek a law of evolution in the form
\begin{equation}
\partial_t \phi=-M\mathrm{grad}_0 E^\mathrm{bulk}(\phi), \label{grad1}
\end{equation}
for some constant mobility $M>0$. To specify the meaning of the constraint gradient denoted by  $\mathrm{grad}_0$, we introduce some function spaces (cf. \cite{Fife10,Kim14}). For every $f\in H^1(\Omega)'$, we denote by $\overline{f}$ its  generalized mean value over $\Omega$ such that
$\overline{f}=|\Omega|^{-1}\langle f,1\rangle_{(H^1)',H^1}$; if $f\in L^1(\Omega)$, then its mean is simply given by $\overline{f}=|\Omega|^{-1}\int_\Omega f(x) \,\d x$. Consider the realization of the minus Laplacian with homogeneous Neumann boundary condition
$\mathcal{A}_N\in \mathcal{L}(H^1(\Omega),H^1(\Omega)')$ defined by
$$ \langle \mathcal{A}_N u,v\rangle_{(H^1)',H^1} := \int_\Omega \nabla u\cdot \nabla v \, \d x,\quad \forall\,u,v\in H^1(\Omega).
$$
For the linear spaces
$$
\dot{H}^1(\Omega)=\{ u \in H^1(\Omega):\ \overline{u}=0\}, \quad
\dot{H}^{-1}(\Omega)= \{ u \in H^1(\Omega)':\ \overline{u}=0 \},
$$
we can check that the operator $\mathcal{A}_N$ is self-adjoint, positively defined on $\dot{H}^1(\Omega)$ and
the restriction of $\mathcal{A}_N$ from $\dot{H}^1(\Omega)$ onto $\dot{H}^{-1}(\Omega)$ is an isomorphism (an easy consequence from the Poincar\'{e}--Wirtinger inequality and the Lax--Milgram theorem). Denote its inverse by $\mathcal{N} =\mathcal{A}_N^{-1}: \dot{H}^{-1}(\Omega) \to \dot{H}^{1}(\Omega)$. Then for every $f\in \dot{H}^{-1}(\Omega)$, $u= \mathcal{N} f \in \dot{H}^{1}(\Omega)$ is the unique weak solution of the Neumann problem
$$
\begin{cases}
-\Delta u=f, \quad \text{in} \ \Omega,\\
\partial_{\mathbf{n}} u=0, \quad \ \  \text{on}\ \partial \Omega.
\end{cases}
$$
Moreover, it holds
\begin{align}
( g, f)_{\dot{H}^{-1}}=
\langle  g, \mathcal{N} f\rangle_{\dot{H}^{-1},\dot{H}^{1}}
=\langle f, \mathcal{N} g\rangle_{\dot{H}^{-1},\dot{H}^{1}} = \int_{\Omega} \nabla(\mathcal{N} g)
\cdot \nabla (\mathcal{N} f) \, \d x, \quad \forall \, g,f \in \dot{H}^{-1}.
\notag
\end{align}
Let
$$
H^4_N(\Omega)=\big\{ \phi\in H^4(\Omega)\ |\ \partial_\mathbf{n}\phi = \partial_\mathbf{n}\Delta \phi=0\ \text{on}\ \partial \Omega\big\}.
$$
We remark that the second boundary condition for $\Delta \phi$ is actually equivalent to $\partial_\mathbf{n}\mu=0$ on $\partial\Omega$ thanks to \eqref{bdN1}.
Let $\phi\in H^4_N(\Omega)$. Then for any smooth function $f$ satisfying $\overline{f}=0$ and $\partial_\mathbf{n} f=0$ on $\partial\Omega$, we use the expression $f=-\Delta(\mathcal{N} f)$ and integration by parts to obtain (in the sense of Gateaux derivative, cf. \cite[Section 3]{Kim14})
\begin{align}
\mathrm{grad}_0 E^\mathrm{bulk}(\phi)[f]
&= \left. \frac{\d}{\d r}E^\mathrm{bulk}(\phi+rf) \right|_{r=0} \notag \\
&= \int_\Omega (\epsilon^2 \nabla \phi\cdot \nabla f +F' (\phi) f) \,\d x \notag \\
&= \int_\Omega (-\epsilon^2 \Delta \phi +F' (\phi)) f \,\d x \notag \\
&= -\int_\Omega (-\epsilon^2 \Delta \phi +F' (\phi)) \Delta (\mathcal{N}f) \,\d x \notag \\
&= \int_\Omega \nabla (-\epsilon^2 \Delta \phi +F' (\phi))\cdot  \nabla (\mathcal{N}f) \,\d x\notag \\
&= \big(-\Delta (-\epsilon^2 \Delta \phi +F' (\phi)), f\big)_{\dot{H}^{-1}}.\notag
\end{align}
Hence, we can identify in $\dot{H}^{-1}(\Omega)$ that
\begin{align}
\mathrm{grad}_0 E^\mathrm{bulk}(\phi)=-\Delta (-\epsilon^2 \Delta \phi +F' (\phi)), \label{grad2}
\end{align}
and specify its domain as $H^4_N(\Omega)$. From \eqref{grad1} and \eqref{grad2} we obtain
$$\partial_t\phi=M\Delta (-\epsilon^2 \Delta \phi +F' (\phi)),$$
 which is exactly the Cahn--Hilliard equation \eqref{claCH} with the mobility $M$ being a positive constant and subject to the boundary conditions \eqref{bdN1}--\eqref{bdN2}.

 The above simple argument indicates that in case of a constant mobility, the Cahn--Hilliard equation can be regarded as a constrained gradient flow in the Hilbert space $\dot{H}^{-1}$. Concerning the more general case with concentration dependent mobilities, it was shown in \cite{Otto01} that for the linear mobility $M(s) = s$, the Cahn--Hilliard equation generates a gradient flow in the space endowed with a non-Hilbertian metric, i.e., the $L^2$-Wasserstein distance (note that the solution considered in \cite{Otto01} is nonnegative in the framework therein). Moreover, in \cite{LMS12}, the authors considered the Cahn--Hilliard type equation with general nonlinear mobilities (e.g., nonnegative concave functions) as a gradient flow in certain weighted-Wasserstein metric spaces and proved the existence of weak solutions by the variational minimizing movement approach.\medskip

\textbf{Derivation by the Energetic Variational Approach}.
We note that the previous derivations, though performed along different procedures, aim to derive a set of partial differential equations that mainly fulfill two physical constraints: the mass conservation (expressed either in a local or a global form) and the energy dissipation. The constitute relation \eqref{fluxJ} is postulated to guarantee the energy dissipation \eqref{BELz}, which may not be the unique choice for this purpose. There are some other evolution equations turning out to have similar properties, for instance, the conservative Allen--Cahn equation:
$$
\partial_t \phi = \epsilon^2 \Delta \phi - F'(\phi) - \frac{1}{|\Omega|}\int_\Omega (\epsilon^2\Delta \phi - F'(\phi))\,\d x,\quad \text{in}\ \Omega\times(0,T),
$$
subject to the homogeneous Neumann boundary condition $\partial_\mathbf{n}\phi=0$ on $\partial \Omega\times (0,T)$. In this case, a direct calculation yields that \eqref{massz} holds and the free energy $E^\mathrm{bulk}(\phi(t))$ is still non-increasing in time:
\begin{align}
\frac{\d}{\d t} E^\mathrm{bulk}(\phi(t))
& =\int_\Omega \left(\epsilon^2 \nabla \phi \cdot \nabla \partial_t \phi
+ F'(\phi) \partial_t \phi \right)\,\d x
\notag
\\
& =\int_\Omega \left(-\epsilon^2 \Delta\phi + F'(\phi)\right) \partial_t \phi \,\d x
\notag
\\
&=\int_\Omega \big(-\partial_t \phi-\overline{ \epsilon^2\Delta \phi - F'(\phi)}\big) \partial_t \phi\, \d x
 \notag
\\
&=-\int_\Omega |\partial_t \phi |^2\,\d x -\overline{ \epsilon^2\Delta \phi - F'(\phi)} \int_{\Omega} \partial_t \phi \,\d x\notag
\\
&= -\int_\Omega |\partial_t \phi |^2 \, \d x,\quad \forall\, t\in (0,T).
\label{BELAC}
\end{align}

Below we provide an alternative derivation of the Cahn--Hilliard equation \eqref{claCH} via the energetic variational approach, which combines the least action principle and Onsager's principle of maximum energy dissipation in continuum mechanics \cite{On31,On31a,Ray}. Within this general framework, one can easily include different physical processes by choosing specific form of the free energy as well as dissipation for the system. Then based on suitable kinematic and energetic assumptions,
the corresponding partial differential equations will be uniquely determined by force balance relations. This approach has been successfully applied to derive complex hydrodynamic systems in fluid dynamics, liquid crystals, electrokinetics, visco-elasticity, reaction-diffusion process and so on, we refer to \cite{DLRW, EHL10, HKL10, KLG, LS03, LW, SV, WLLE20, WXL13, XSL14} and the references cited therein.

In $\Omega$, $\phi$ is assumed to be a locally conserved quantity that satisfies the continuity equation
\begin{equation}
\partial_t \phi +\nabla\cdot(\phi \mathbf{u})=0,
\quad (x, t)\in \Omega\times (0,T).
\label{con1z}
\end{equation}
Here, the mass flux is given by the advection flux $\mathbf{J}=\phi \mathbf{u}$, where $\mathbf{u}: \Omega \to \mathbb{R}^d$ stands for the microscopic effective velocity (e.g., due to certain diffusion process). Concerning the boundary conditions, we again assume the no-flux boundary condition $\mathbf{J}\cdot \mathbf{n}=0$ on $\partial \Omega$, which can be guaranteed by
\begin{equation}
\mathbf{u}\cdot \mathbf{n}=0,
\quad (x, t)\in \partial \Omega \times (0,T).
\label{bd1z}
\end{equation}

For an isothermal closed system, we assume that the evolution of a binary mixture satisfies the following energy dissipation law, which is a natural consequence of the first and second law of thermodynamics:
\begin{align}
 \frac{\d}{\d t}E^{\mathrm{bulk}}(t)=-\mathcal{D}^{\mathrm{bulk}}(t),\quad t\in (0,T),
 \label{ABELz}
\end{align}
where the bulk free energy $E^{\mathrm{bulk}}$ takes the form in \eqref{Ebulk} and the rate of energy dissipation $\mathcal{D}^{\mathrm{bulk}}$ is chosen as
\begin{align}
\mathcal{D}^{\mathrm{bulk}}(t)= \int_\Omega M^{-1}\phi^2 |\mathbf{u}|^2 \, \d x,
\label{diss1z}
\end{align}
with $M>0$ being the mobility.

It remains to determine the microscopic velocity $\mathbf{u}$ in \eqref{con1z}. Let $\Omega_0^X, \Omega_t^x\subset\mathbb{R}^d$ be bounded domains with smooth boundaries $\Gamma_0^X, \Gamma_t^x$, respectively.
We introduce the flow map  $x(X,t): \Omega_0^X\to \Omega_t^x$, which is defined as a solution to the system of ordinary differential equations
\begin{equation}
\begin{cases}
\displaystyle{\frac{\d}{\d t}}x(X,t)=\mathbf{u}(x(X,t),t),\quad t>0,\\
x(X,0)=X,
\end{cases}
\label{flowmap}
\end{equation}
where $X=(X_1,...,X_d)^\mathrm{tr}\in \Omega_0^X$, $x=(x_1,...,x_d)^\mathrm{tr}\in \Omega_t^x$, and  $\mathbf{u}(x,t)\in \mathbb{R}^d$ is a (smooth) velocity field. The coordinate system $X$ is called the Lagrangian coordinate system and refers
to $\Omega_0^X$ that we call the reference configuration; the coordinate system
$x$ is called the Eulerian coordinate system and refers to $\Omega_x^t$
that we call the deformed configuration. We shall denote $\nabla_x$ the gradient in $\Omega$ under the Eulerian coordinate system. Introduce the action functional
\begin{align}
\mathcal{A}^{\mathrm{bulk}}(x(X,t))  = -\int_0^T  E^{\mathrm{bulk}}(\phi(t))\,\d t.
\label{ac1}
\end{align}
Then applying the least action principle, which states that the dynamics of a Hamiltonian
system is determined by a critical point of the action functional with respect to the trajectory (in Lagrangian coordinates), we get
\begin{align}
\delta_{x} \mathcal{A}^{\mathrm{bulk}}
&  = -\int_0^T\!\int_{\Omega_t^x} (\phi \nabla_x \mu) \cdot \delta x\, \d x \d t,\quad \text{with}\ \ \mu=-\epsilon^2\Delta \phi +F'(\phi).
\label{LAP1}
\end{align}
In the above computation, we have assumed a trivial boundary dynamics for $\phi$ such that the boundary condition \eqref{bdN2} holds. The relation \eqref{LAP1} yields the conservative force (in the strong form, cf. \cite[Section 2.2.1]{SV})
$$
\mathbf{f}_{\text{conv}}=-\phi \nabla \mu.
$$
On the other hand, from Onsager's maximum dissipation principle for a dissipative system, we can calculate the dissipative force by taking variation of the Rayleigh dissipation functional $\mathcal{R}=\frac12\mathcal{D}^{\mathrm{bulk}}$ (recall \eqref{diss1z}) with respect to the rate function $\mathbf{u}$, that is,
\begin{align}
\delta_{\mathbf{u}} \left(\frac12 \mathcal{D}^{\mathrm{bulk}}\right)
&= \frac12\left.\frac{\d}{\d r} \right|_{r=0}  \mathcal{D}^{\mathrm{bulk}}(\mathbf{u}+r \tilde{\mathbf{u}})= \int_{\Omega} M^{-1} \phi^2 \mathbf{u}\cdot \delta \mathbf{u}\, \d x.\nonumber
\end{align}
This gives the generalized dissipative force (again written in the strong form, cf. \cite[Section 2.2.1]{SV}):
$$
\mathbf{f}_{\text{diss}}=-M^{-1}\phi^2 \mathbf{u}.
$$
By the force balance relation, i.e., Newton's second law  $\mathbf{f}_{\text{inertial}}+\mathbf{f}_{\text{conv}}+\mathbf{f}_{\text{diss}}=0$ (recalling that here we have $\mathbf{f}_{\text{inertial}}=0$ because the kinetic energy is neglected), we obtain
\begin{equation*}
\phi \nabla \mu +M^{-1}\phi^2 \mathbf{u}=0, \quad \text{in}\ \Omega\times(0,T),
\end{equation*}
where $\mu =-\epsilon^2 \Delta \phi +F'(\phi)$. Finally, solving $\mathbf{u}$ from the above algebraic equation and inserting it back into \eqref{con1z}, we arrive at the Cahn--Hilliard equation \eqref{claCH}, subject to the boundary conditions \eqref{bdN1}--\eqref{bdN2}.
The above derivation via the energetic variational approach reveals that the Cahn--Hilliard equation together with the classical boundary conditions naturally fulfills three important physical constraints: conservation of mass, dissipation of energy and, in addition, force balance.

\subsection{Well-posedness}
The Cahn--Hilliard equation \eqref{claCH} subject to the initial and boundary conditions \eqref{ini0}--\eqref{bdN2} has been extensively studied in the literature. A rather complete picture about the existence, the uniqueness, the regularity and the long-time  behavior of solutions has been obtained. We refer the reader to the recent book \cite{Mi19} and the references cited therein.\medskip

\textbf{The case with regular potential}. Concerning problem \eqref{claCH}--\eqref{bdN2} with a constant mobility (e.g., $M=1$ without loss of generality) and a regular potential $F$ (e.g., a fourth order polynomial like \eqref{regpo}),
 existence of global weak solutions with an initial datum $\phi_0\in H^1(\Omega)$ can be easily obtained by making use of a Faedo--Galerkin method, thanks to the Lyapunov structure \eqref{BELz} (see \cite{NS85,E89}). Existence and uniqueness of more regular solutions with initial datum $\phi_0\in H^2_N(\Omega)=\{\phi\in H^2(\Omega)\,|\, \partial_\mathbf{n}\phi=0\ \text{on}\ \Omega\}$ was proved in \cite{EZ86}. More precisely, they consider
$$
\mu=-\epsilon^2 \Delta \phi + F'(\phi),\quad \text{with}\ \ F'(\phi)=\gamma_2 \phi^3+\gamma_1\phi^2-\phi,
$$
where $\gamma_1,\gamma_2$ are constant parameters. It was shown that the sign of $\gamma_2$ (coefficient of the leading order term) plays a crucial role in the study of existence of global solutions. If $\gamma_2>0$, for any initial datum there is a unique global solution, while for $\gamma_2<0$, the solution must blow up in a finite time for large initial data. The proofs therein rely on the Picard iteration scheme (for local existence and uniqueness) and the energy method (for global existence and finite time blow-up). Besides, they studied a finite element Galerkin approximation of the initial boundary value problem and obtained existence results as well as optimal order error bounds. We also refer to \cite{NST89}, in which the authors considered a more general polynomial form of $F'(s)=\sum_{j=1}^{2p-1}a_js^j$ ($p\in\mathbb{N}, p\geq 2$) and boundary conditions of either Neumann or periodic type. The form and growth condition on the regular potential $F$ can be relaxed when studying the regularity of solutions, see for instance, \cite[Chapter 3, Section 3.4]{Mi19}: $F\in C^3(\mathbb{R})$, $F(0)=F'(0)=0$ and
\begin{align*}
& F''(s)\geq -c_0,\quad c_0\geq 0, \quad \forall\, s\in \mathbb{R},\\
& F'(s)s\geq c_1 F(s)-c_2,\quad F(s)\geq -c_3,\quad c_1>0,\ c_2,c_3\geq 0,\quad \forall\, s\in\mathbb{R},\\
& |F'(s)|\leq \varepsilon F(s)+c_\varepsilon,\quad \forall\, \varepsilon>0, \ s\in\mathbb{R}.
\end{align*}
It is worth noting that when the mobility is a constant, \eqref{claCH} is a fourth-order semilinear parabolic equation. In general, a fourth-order parabolic equation does not enjoy the maximum principle for second-order parabolic equations. Concerning \eqref{claCH}, its solution $\phi$ needs not to stay in the physically relevant interval $[-1,1]$ along evolution, even if the initial datum $\phi_0$ satisfies this property. We refer to \cite{P11} as well as \cite[Chapter 4, Remark 4.10]{Mi19} for a simple counterexample regarding the Cahn--Hilliard equation with $M=1$ and $F'(\phi)=\phi^3-\phi$ in one space dimension.

Next, we would like to mention the viscous Cahn--Hilliard equation
\begin{equation}
\left\{
\begin{aligned}
&\partial_t\phi=\Delta \mu,\quad &\text{in}\ \Omega \times (0,T),\\
&\mu=\alpha \partial_t \phi-\epsilon^2\Delta \phi + F'(\phi),\quad &\text{in}\ \Omega \times (0,T),
\end{aligned}
\right.
\label{clavCH}
\end{equation}
with the viscous parameter $\alpha >0$. Equation \eqref{clavCH} was introduced in \cite{NC88} to include certain viscoelastic relaxation effects in the phase separation process, which was neglected in \cite{CH58}. The viscous term $\alpha \partial_t \phi$ is also related to the influence of certain internal microforces in the mixture (see \cite{Gur96}). From the energetic point of view, it yields some additional dissipation in the system, for instance, any smooth solution to the initial boundary value problem \eqref{clavCH}, \eqref{ini0}--\eqref{bdN2} satisfies
\begin{align}
\frac{\d}{\d t} E^\mathrm{bulk}(\phi(t))
= -\int_\Omega |\nabla \mu|^2 \, \d x-\alpha \int_\Omega |\partial_t \phi|^2\,\d x,\quad \forall\, t\in (0,T).
\label{BELzv}
\end{align}
Thus, the viscous Cahn--Hilliard equation \eqref{clavCH} can be regarded as a regularized version of the original Cahn--Hilliard equation. When the potential $F$ is regular, i.e., a general polynomial, we refer to \cite{BEG95,ES96} for an extensive study on computations and mathematical analysis for the equation \eqref{clavCH} subject to the homogeneous Dirichlet (or Neumann) boundary conditions as well as the initial condition \eqref{ini0}. The results in \cite{BEG95,ES96} showed that in a suitable sense the viscous Cahn--Hilliard equation can actually be viewed as an interpolation between the Cahn--Hilliard equation \eqref{claCH} (with $M=1$) and the Allen--Cahn equation \cite{AC79} for grain-boundary migration, namely, $\partial_t \phi=\epsilon^2\Delta \phi - F'(\phi)$. \medskip

\textbf{The case with degenerate mobility}. Next, we consider the situation that the mobility $M$ is concentration dependent, and even may be degenerate at some values of $\phi$. This case appeared in the original derivation of the Cahn--Hilliard equation \cite{CH61}. A thermodynamically reasonable choice is $M(s)=1-s^2$ (see \cite{Hi61} and also \cite{NS84,CT94,NC08}). Concerning the mathematical analysis, the author of \cite{Yin92} studied a general Cahn--Hilliard type equation in one space dimension (i.e., on the interval $(0,1)$) with $F'$ being continuous and $M$ being H\"{o}lder continuous  such that
$$
M(0)=M(1)=0, \quad M(s)\geq 0, \ \ \forall\, s\in (0,1).
$$
Then for any initial datum $\phi_0\in H^3_0(0,1)=\{\phi\in H^3(0,1)\,|\, Du=0\ \text{at}\ x=0,1\}$ satisfying $0\leq \phi_0\leq 1$, he proved the existence of a global generalized solution $\phi$ with uniform bounds in $L^\infty(0,T;H^1(0,1))$ and local $L^2$ estimates on $D^3\phi$ when the equation does not degenerate. In particular, thanks to the degeneracy of the mobility, he obtained nonnegativity of the solution $0\leq \phi(x,t)\leq 1$ (note that in \cite{Yin92}, $\phi$ stands for the concentration of one component, not the difference), see \cite[Theorem 3.2]{Yin92}. The result of \cite{Yin92} implies that the degenerate mobility turns out to be a sufficient condition for the existence of weak solutions with the physical property $0\leq \phi(x,t)\leq 1$ for $t\geq 0$, as long as $0\leq \phi_0(x)\leq 1$.

Later on, the Cahn--Hilliard equation \eqref{claCH} with a degenerate concentration dependent mobility as well as a singular potential was analyzed in \cite{EG96}. Assume that (i) $F(s)=F_1(s)+F_2(s)$, with $\|F_2(s)\|_{C^2[-1,1]}\leq C$ and $F_1:(-1,1)\to \mathbb{R}$ being convex, satisfying $F''_1(s)=(1-u^2)^{-m}\widetilde{F}(s)$ for $m\geq 1$ and a $C^1$-function $\widetilde{F}:[-1,1]\to \mathbb{R}^+\cup\{0\}$; (ii) $M(s)=(1-s^2)^m\widetilde{M}(s)$ with a $C^1$-function $\widetilde{M}:[-1,1]\to \mathbb{R}$ satisfying $0<m_0\leq \widetilde{M}(s)\leq M_0$ for $s\in [-1,1]$, and $M(s)$ is extended to $\mathbb{R}$ by setting $M(s)=0$ for $|s|>1$. It is easy to check that these general assumptions cover the physically relevant case with the logarithmic potential \eqref{logpo} and $M(s)=1-s^2$, $s\in[-1,1]$. Let either $\partial \Omega\in C^{1,1}$ or $\Omega$ be convex in $\mathbb{R}^d$ ($d\geq 1$). Then for any initial datum $\phi_0\in H^1(\Omega)$ with $|\phi_0|\leq 1$ a.e. in $\Omega$ and $\int_\Omega F(\phi_0)+ \Phi(\phi_0)\, \mathrm{d}x<+\infty$, problem \eqref{claCH}--\eqref{bdN2} admits a global weak solution on an arbitrary time interval $[0, T]$, see \cite[Theorem 1]{EG96}. Here, the function $\Phi:(-1,1)\to \mathbb{R}^+\cup \{0\}$ is determined by $\Phi''(s)=M(s)^{-1}$, $\Phi(0)=\Phi'(0)=0$. The proof relies on some suitable approximation of the degenerate mobility and regularization of the singular potential such that (see \cite[Section 3]{EG96})
\begin{equation*}
M_\varepsilon(s):=
\begin{cases}
M(-1+\varepsilon),\quad \text{for}\ s\leq -1+\varepsilon,\\
M(s),\qquad\quad \ \, \text{for}\ |s|\leq 1-\varepsilon,\\
M(1-\varepsilon),\qquad\! \text{for}\ s\geq 1-\varepsilon,
\end{cases}
\end{equation*}
$F_\varepsilon=F^1_\varepsilon+F_2$ with
\begin{equation*}
(F^1_\varepsilon)''(s):=
\begin{cases}
F_1''(-1+\varepsilon),\quad \text{for}\ s\leq -1+\varepsilon,\\
F_1''(s),\qquad\quad \ \,  \text{for}\ |s|\leq 1-\varepsilon,\\
F_1''(1-\varepsilon),\qquad\! \text{for}\ s\geq 1-\varepsilon,
\end{cases}
\quad
F^1_\varepsilon(0)=F_1(0), \quad (F^1_\varepsilon)'(0)=F_1'(0).
\end{equation*}
and $F_2$ being extended to be a function on $\mathbb{R}$ such that $\|F_2\|_{C^2(\mathbb{R})}\leq C$.
We note that the special structure that $M(s)F''(s)$ is bounded plays a crucial role in the proof of existence in \cite{EG96}.

Recently, the authors of \cite{DD16} considered the Gibbs--Thomson effect in the phase separation process. For the Cahn--Hilliard equation \eqref{claCH} in $\Omega=[0,2\pi]^d$ subject to periodic boundary conditions, assuming that the double-well potential $F$ is smooth at its minima $\pm1$ (cf. \eqref{regpo}) and the mobility is of the form $M(s)=|1-s^2|^m$ \textit{without} zero extension outside $[-1,1]$ ($m>0$ when $d=1,2$, and $m\in (0,2/(d-2))$ when $d\geq 3$), they proved that for any initial datum $\phi_0\in H^1(\Omega)$, the resulting problem admits a global weak solution on an arbitrary time interval $[0, T]$, see \cite[Theorem 2]{DD16}. The proof is again based on some regularization of the degenerate mobility and then passing to the limit in the non-degenerate approximating problem. One difference from the previous results (cf. those in \cite{EG96}) is that for the solution $\phi$, even if its initial value lies in $[-1,1]$, it may not remain inside $[-1,1]$ as long as the interface has nonzero mean curvature, which accommodates the physical Gibbs--Thomson effect. Nevertheless, it is worth mentioning that in all studies above, the uniqueness of weak solutions (in various formulations) to the Cahn--Hilliard equation with degenerate mobility remains an open problem. \medskip

\textbf{The case with singular potential}. A lot of attentions have been paid to problem \eqref{claCH}--\eqref{bdN2} with a constant mobility $M>0$ and a singular potential $F$ including the typical ones \eqref{logpo} and \eqref{obspo}, see \cite{CMZ11,Mi17,Mi19} for some surveys on this topic.

To study the existence of solutions, one basic strategy is to regularize the singular potential in a suitable manner, prove the existence of solutions to the regularized problem, derive uniform estimates with respect to the approximating parameter and then pass to the limit to extract a convergent subsequence.

In \cite{BE91}, the authors studied problem \eqref{claCH}--\eqref{bdN2} in a bounded smooth domain $\Omega\subset \mathbb{R}^d$ ($d=1,2,3$) with the double obstacle potential \eqref{obspo}. They proved that for any initial datum $\phi_0\in H^1(\Omega)$, $|\phi_0(x)|\leq 1$ a.e. in $\Omega$ with its spatial average in $(-1,1)$, problem \eqref{claCH}--\eqref{bdN2} admits a unique global weak solution $\phi$ satisfying $-1\leq \phi(x,t)\leq 1$ a.e. in $\Omega \times (0,T)$, see \cite[Theorem 2.2]{BE91}. Some regularity results on the weak solutions were also obtained in \cite[Section 2.2]{BE91}. The proof relies on the following $C^2$-regularization for the singular potential (see \cite[(2.9)]{BE91}): for $\varepsilon \in (0,1)$,
\begin{align*}
F_\varepsilon(s)=
\begin{cases}
\displaystyle{\frac{1}{2\varepsilon}\Big(s-\Big(1+\frac{\varepsilon}{2}\Big)\Big)^2+\frac12(1-s^2) +\frac{\varepsilon}{24},\quad\text{for}\ s\geq 1+\varepsilon,}\smallskip \\
\displaystyle{\frac{1}{6\varepsilon^2}(s-1)^3 +\frac12(1-s^2),\qquad\qquad\qquad\ \, \text{for}\ 1<s< 1+\varepsilon,}\smallskip\\
\displaystyle{\frac12(1-s^2),\qquad\qquad\qquad\qquad\qquad\qquad\ \ \,\text{for}\ |s|\leq 1,}\smallskip\\
\displaystyle{-\frac{1}{6\varepsilon^2}(s+1)^3+\frac12(1-s^2),\quad \qquad\qquad\ \ \, \text{for}\ -1-\varepsilon<s<-1,}\smallskip\\
\displaystyle{\frac{1}{2\varepsilon}\Big(s+\Big(1+\frac{\varepsilon}{2}\Big)\Big)^2 +\frac12(1-s^2) +\frac{\varepsilon}{24},\quad\,\text{for}\ s\leq -1-\varepsilon.}
\end{cases}
\end{align*}

When the logarithmic potential \eqref{logpo} is concerned, several different methods have been developed in the literature. The authors of \cite{EL91} studied a general multi-component problem in a bounded smooth domain $\Omega\subset \mathbb{R}^d$ ($d=1,2,3$) subject to homogeneous Neumann boundary conditions. Their results can easily apply to the two-component problem, that is, for potentials like $F(s)=\theta [s\ln s+ (1-s)\ln(1-s)]-(a_0s^2+a_1s+a_2)$, $s\in (0,1)$, if any initial datum $\phi_0\in H^1(\Omega)$, $0\leq \phi_0(x)\leq 1$ a.e. in $\Omega$ with its spatial average in $(0,1)$, problem \eqref{claCH}--\eqref{bdN2} admits a unique global weak solution $\phi$ satisfying $0\leq \phi(x,t)\leq 1$ a.e. in $\Omega \times (0,T)$, see \cite[Theorem 1]{EL91} for a general statement of the result. Indeed, the arguments in \cite[Section 3]{EL91} also implied that $0< \phi(x,t)< 1$ a.e. in $\Omega \times (0,T)$, i.e., the set of singular points for $F'(s)$ has zero measure. Besides, in \cite[Theorem 2]{EL91} they justified the ``deep-quench limit" problem studied in \cite{BE91} by letting $\theta\to 0^+$. The key idea of the proof in \cite{EL91} is to regularize the function $G(s)=s\ln s$ by
\begin{align*}
G_\varepsilon (s)=
\begin{cases}
s\ln s,\qquad \quad \qquad\ \ \text{for}\ s\geq \varepsilon,\\
\displaystyle{\frac{s^2}{2\varepsilon} +s\ln \varepsilon -\frac{\varepsilon}{2},}\quad\ \,  \text{for}\ s<\varepsilon,
\end{cases}
\end{align*}
for some $\varepsilon\in (0,1)$. After solving the regularized problem, they derived estimates that are uniform in $\varepsilon$ and then passed to the limit as $\varepsilon\to 0^+$ to extract a convergent subsequence, whose limit is a global weak solution to the original problem.

Later in \cite{DD95}, the authors introduced a different approximation of the logarithmic potential \eqref{logpo} in terms of polynomials
$$
F_n(s)=\theta \sum_{k=0}^n\frac{s^{2k+2}}{(2k+1)(2k+2)}-\frac{\theta_c}{2}s^2,\quad \forall\,s\in (-1,1), \ n\in \mathbb{N}.
$$
Then for any initial datum $\phi_0\in L^2(\Omega)$ (or $H^1(\Omega)$), $\|\phi_0\|_{L^\infty(\Omega)}\leq 1$ with its spatial average in $(-1,1)$, they prove that problem \eqref{claCH}--\eqref{bdN2} (or in a periodic setting) admits a unique global weak solution satisfying $\|\phi(t)\|_{L^\infty(\Omega)}\leq 1$ for $t\geq 0$ and for $t>0$, the set $\{x\in \Omega\,|\,|\phi(x,t)|=1\}$ has measure zero, see \cite[Theorem 2.2]{DD95}. Moreover, they showed the continuity property of weak solutions, which yields the existence of $C_0$-semigroups defined by $S(t):\phi_0\to \phi(t)$ in suitable phase spaces.

In \cite{Mi17}, the author considered the following equation
\begin{equation}
\left\{
\begin{aligned}
&\partial_t\phi+\chi \phi =\Delta \mu,\quad &\text{in}\ \Omega \times (0,T),\\
&\mu=-\epsilon^2\Delta \phi + F'(\phi),\quad &\text{in}\ \Omega \times (0,T),
\end{aligned}
\right.
\label{claCHO}
\end{equation}
with some $\chi\geq 0$ and a general class of singular nonlinearities $F'\in C^1(-1,1)$ satisfying $F(0)=F'(0)=0$ and
\begin{equation}
\lim_{s\to\pm 1}F'(s)=\pm \infty,\quad \lim_{s\to \pm 1} F''(s)=+\infty.
\label{sing1}
\end{equation}
In particular, the logarithmic potential \eqref{logpo} fulfills the above assumptions. By adapting the methods in \cite{EL91,EG96}, he made use of the following approximation $F'_\varepsilon\in C^1(\mathbb{R})$: for $\varepsilon\in (0,1)$,
\begin{equation*}
F'_\varepsilon(s)=
\begin{cases}
F'(-1+\varepsilon)+F''(-1+\varepsilon)(s+1-\varepsilon),\quad \text{for}\ s<-1+\varepsilon,\\
F'(s),\qquad \qquad\qquad\qquad \qquad\qquad\qquad\  \text{for}\ |s|\leq 1-\varepsilon,\\
 F'(1-\varepsilon)+F''(1-\varepsilon)(s-1+\varepsilon),\qquad\ \ \text{for}\ s>1-\varepsilon.
\end{cases}
\end{equation*}
Let $\Omega\subset \mathbb{R}^d$, $d=1,2,3$ be a smooth bounded domain. He proved that for any initial datum $\phi_0\in H^1(\Omega)$, $|\phi_0(x)|< 1$ a.e. in $\Omega$ with its spatial average in $(-1,1)$, the initial boundary value problem \eqref{claCHO}, \eqref{ini0}--\eqref{bdN2} admits a unique global weak solution $\phi$ satisfying $-1< \phi(x,t)< 1$ a.e. in $\Omega \times (0,T)$, see \cite[Theorem 2.6, Remark 2.7, Remark 2.8]{Mi17}. We note that when $\chi>0$, the equation \eqref{claCHO} is usually referred to as the Cahn--Hilliard--Oono equation, which was introduced to model long-ranged (nonlocal) interactions \cite{OP87}. For further details on its mathematical analysis, see e.g., \cite{Mi11} (with the regular potential \eqref{regpo}) and \cite{GGM17} (with singular potentials including \eqref{logpo}). In \cite{Mi11}, it was also shown that, in certain sense (i.e., the existence of a robust family of exponential attractors as $\chi\to0^+$), the dynamics of the Cahn--Hilliard equation is ``close" to that of the Cahn--Hilliard--Oono equation, for $\chi>0$ being small.

Besides the approximating arguments described above, there are some other approaches to handle the Cahn--Hilliard equation with singular potentials. We first mention the work \cite{KNI95}, which treated problem \eqref{claCH}--\eqref{bdN2} (again with a constant mobility) in a quite general setting. For a wide class of non-smooth potentials including both \eqref{logpo} and \eqref{obspo}, thanks to the decomposition
$F(s)=\widehat{\beta}(s)+\widehat{g}(s)$, where $\widehat{\beta}$ is a proper, lower semi-continuous and convex function and $\widehat{g}$ is a $C^1$ function with Lipschitz continuous derivative $g=(\widehat{g})'$, the Cahn--Hilliard equation \eqref{claCH} can be written in a generalized form involving multi-valued mappings via subdifferential operators. In this framework, the authors first regularized the Cahn--Hilliard equation by the viscous one (cf. \eqref{clavCH}) written in the following form
\begin{align*}
\begin{cases}
(\alpha I+\mathcal{N})\partial_t\phi=\epsilon^2\Delta \phi-F'(\phi)+|\Omega|^{-1}\int_\Omega F'(x)\,\d x,\quad \text{in}\ \Omega\times (0,T),  \\
\partial_\mathbf{n}\phi=0,\qquad \qquad \qquad \qquad \qquad \qquad \qquad \qquad \qquad \quad  \text{on}\ \partial \Omega\times (0,T),
\end{cases}
\end{align*}
where $\alpha>0$ and $\mathcal{N}$ denotes the inverse of minus Laplacian subject to the homogeneous Neumann boundary condition on the space $L^2_0(\Omega)=\{\phi\in L^2(\Omega)\ |\  \int_\Omega \phi(x)\,\d x=0\}$. Then they solved the viscous problem (indeed in a more general form) using the abstract result in \cite{KP86}, whose proof was based on Yosida's approximation of monotone operators. After deriving uniform estimates of approximate solutions with respect to the viscous parameter $\alpha$, they obtained the existence of a global weak solution to the original problem \eqref{claCH}--\eqref{bdN2} by finding a convergent subsequence as $\alpha\to 0^+$ (see \cite[Theorem 6.1]{KNI95} for a complete statement including also results on uniqueness and long-time behavior). On the other hand, with the help of monotone operator methods \cite{Sh97}, the authors of \cite{AW07} gave a different proof for the existence and uniqueness of global weak solutions to problem \eqref{claCH}--\eqref{bdN2} with logarithmic type potentials (see \cite[Theorem 1.2]{AW07}). After decomposing the nonlinear singular term $F'$ in  \eqref{claCH} into a monotone operator plus a globally Lipschitz continuous one, they achieved the conclusion in a direct manner by solving an abstract Cauchy problem for a suitable Lipschitz perturbation of a monotone operator, which is the subgradient of the convex part of the energy $E^{\mathrm{bulk}}(\phi)$ (see \cite[Theorem 3.1]{AW07}). \medskip

\textbf{Regularity and separation from pure states}.
Regularity of solutions to problem \eqref{claCH}--\eqref{bdN2} has been investigated in some of the works mentioned above and the references cited therein. Roughly speaking, the parabolic nature of the Cahn--Hilliard equation (valid when $M>0$) yields some instantaneous regularizing effect of its weak solutions for $t>0$. This fact is well understood for the case with regular potentials, while the case with singular potentials like \eqref{logpo} and \eqref{obspo} turns out to be more tricky. On one hand, the possible singularity of the potential (or its derivatives) at the pure states $\pm 1$ guarantees that the order parameter $\phi$ stays in the physically relevant interval $(-1,1)$ (cf. the case with degenerate mobility). On the other hand, it brings further difficulties to gain higher-order spatial regularity of the solution, see extensive discussions in   \cite{AW07,BE91,CMZ11,DD95,KNI95,Mi17,Mi19,MZ04}.

The case with a logarithmic potential like \eqref{logpo} is of particular interest, because its two minima $\pm\phi_*$, which are the two nonzero solutions to the equation $\ln(1+s)-\ln(1-s)=(2\theta_c/\theta) s$ with $0<\theta<\theta_c$, locate exactly inside the interval $(-1,1)$. From the dissipative nature of the Cahn--Hilliard equation, if the initial value is not a pure state (e.g., $|\Omega|^{-1}|\int_\Omega \phi_0(x)\,\d x|<1$), one may expect in this case the following strict separation property along evolution:
\begin{equation}
\|\phi(t)\|_{L^\infty(\Omega)}\leq 1-\delta, \label{sep}
\end{equation}
for some $\delta\in (0,1)$. The property \eqref{sep}, if it holds, implies that in the conserved phase separation process governed by the Cahn--Hilliard equation, the pure states can never be completely reached. From the mathematical point of view, it plays an important role in the analysis of \eqref{claCH}, since the singular potential \eqref{logpo} can thus be regarded as a smooth, globally Lipschitz function so that further regularity of solutions to \eqref{claCH} can be obtained, see e.g., \cite{AW07,Mi17,Mi19,MZ04}.

In \cite{MZ04}, the authors considered the viscous Cahn--Hilliard equation \eqref{claCHO} when the spatial dimension is less or equal than three. Taking advantage of the viscous term $\alpha \partial_t \phi$ ($\alpha>0$) and using the comparison principle for second-order parabolic equations, they proved that \eqref{sep} holds for all $t > 0$ (i.e., the so-called instantaneous strict separation property) and the separation is uniform when $t\geq \eta$, for an arbitrary but fixed $\eta > 0$ (see \cite[Corollary 3.2]{MZ04}). However, the distance $\delta$ depends on the viscous parameter $\alpha$ so that one cannot pass to the limit as $\alpha\to 0^+$ to draw the same conclusion for the original Cahn--Hilliard equation \eqref{claCH}.

When $\alpha=0$, the spatial dimension indeed plays a crucial role in the argument. In one and two dimensions, the instantaneous strict separation property was proved in \cite[Proposition 7.1, Theorem 7.2]{MZ04} for general singular potentials satisfying \eqref{sing1}. When the dimension is two, an additional assumption was also required therein:
$$|F''(s)|\leq e^{C_1|F'(s)|+C_2}$$
with some positive constants $C_1$, $C_2$, which allows one to apply the Trudinger--Moser type inequality \cite{NSY}. Nevertheless, the above assumptions are fulfilled by the logarithmic potential \eqref{logpo}. Later on, the authors of \cite{GGM17} extended the result to the two dimensional Cahn--Hilliard--Oono equation \eqref{claCHO} by employing an alternative approach, which replies on the analysis of the Neumann problem of an elliptic equation with a singular nonlinearity:
\begin{equation*}
\begin{cases}
-\epsilon^2 \Delta \phi+ \widetilde{F}'(\phi)=f,\quad \text { in } \Omega, \\
\partial_{\mathbf{n}} \phi=0, \qquad \qquad\qquad  \text { on } \partial \Omega,
\end{cases}
\end{equation*}
where the nonlinear function $\widetilde{F}$ is the strictly convex part of $F$  satisfying $\lim_{s\to \pm 1} \widetilde{F}(s)=\pm \infty$ (see \cite[Section 3]{GGM17}, or \cite[Appendix]{GGW18}). We also refer to \cite{GMT19,GGW18,HW21,SW20} for its generalizations to higher order Cahn--Hilliard type equation and some more complicated systems with fluid interaction.

When the spatial dimension is three, the situation is less satisfactory, because the singularity of the logarithmic potential at $\pm 1$ seems not strong enough.
In \cite{MZ04}, the authors obtained the instantaneous strict separation property, under a stronger assumption on the potential such that \eqref{sing1} is satisfied together with $$|F''(s)|\leq C(|F'(s)|^2+1)$$ for some $C>0$  (see \cite{LP18} for some recent improvements). The above assumption is valid, for instance, for a class of nonlinearities like (see \cite[Remark 7.1]{MZ04})
$$
\frac{h(s)}{(1-s^2)^\gamma},\quad \text{where}\ h\in C^1([-1,1]),\ \ h(\pm 1)\neq 0\ \text{and}\ \gamma\geq 1,
$$
but it is not satisfied by \eqref{logpo}, the physically important one. On the other hand, in \cite{AW07}, using the strict separation property of solutions to the stationary problem of the Cahn--Hilliard equation (see \cite[Proposition 6.1]{AW07}), the authors proved that global weak solutions of problem \eqref{claCH}--\eqref{bdN2} with logarithmic type potential will stay eventually away from the pure states $\pm 1$ for sufficiently large time. Besides, the stability result for the Cahn--Hilliard--Hele--Shaw system that was obtained in \cite[Theorem 2.5]{GGW18} implies that (after simply neglecting the fluid interaction), if the initial datum $\phi_0$ is not a pure state and belongs to a suitably small $H^2$-neighborhood of a local minimizer of $E^{\mathrm{bulk}}(\phi)$, then the unique global strong solution to problem \eqref{claCH}--\eqref{bdN2} with logarithmic type potential exists and is strictly separated from $\pm 1$ for all $t\geq 0$.

\subsection{Long-time behavior}

In what follows, we review some results on the long-time behavior of global solutions to problem \eqref{claCH}--\eqref{bdN2} as $t\to +\infty$. Generally speaking, the study of long-time behavior of the global solution to a nonlinear evolution equation can be divided into two categories: the first category is to investigate the long-time behavior of the solution corresponding to a given initial datum, while the second category is to investigate the long-time behavior of a bundle of solutions whose initial data vary in a bounded set (see \cite{Z}). The latter category is related to the study of the associated infinite dimensional dynamical system \cite{Te}, for instance, the existence of (finite dimensional) global attractors, exponential attractors, inertial manifolds, and so on. In this direction, there is a huge literature on problem \eqref{claCH}--\eqref{bdN2} and its variants (for instance, non-isothermal model, hyperbolic relaxation etc), we refer to \cite{CMZ11,DD95,Ga08,GGM17,GPS07,GSSZ09,GSZ09,GSZ10,Mi11,Mi17,MZ04,NS85,NST89,SZ02,Sh07,Te,ZM04,ZM05} and the references cited therein, in particular, we refer to the recent book \cite{Mi19} for more detailed information. In this paper, we confine ourselves to the discussion on the long-time behavior of a single trajectory defined by the global solution to problem \eqref{claCH}--\eqref{bdN2}.

In view of \eqref{BELz}, it is natural to ask whether the global solution to the evolution problem \eqref{claCH}--\eqref{bdN2} will converge to an equilibrium (e.g., local or global minimizers of the energy $E^{\mathrm{bulk}}(\phi)$) as $t\to+\infty$. This property is sometimes also called the uniqueness of asymptotic limit as $t\to +\infty$ for a global solution.
In the one dimensional case, the long-time behavior of global solutions to problem \eqref{claCH}--\eqref{bdN2} with a regular potential like \eqref{regpo} was first analyzed in \cite{Z86}. Based on the well-posedness result obtained in \cite{EZ86} and the Lyapunov structure \eqref{BELz}, the author proved that for any initial datum $\phi_0\in  H^4_N(\Omega)=\{\phi\in H^4(\Omega)\,|\, \partial_\mathbf{n}\phi=\partial_\mathbf{n}\Delta\phi=0\ \text{on}\ \partial \Omega\}$,  the trajectory of solution is relatively compact in $H^2(\Omega)$ and as $t\to+\infty$, $\phi(x,t)$ converges to the $\omega$-limit set of $\phi_0$ denoted by $\omega(\phi_0)$, which is a compact, connected subset of $H^2(\Omega)$ and is positive invariant under the nonlinear semigroup $S(t):\phi_0 \to \phi(t)$ defined by the solution (see \cite[Theorem 2.2]{Z86}). Besides, he proved that each element of $\omega(\phi_0)$ is an equilibrium and in the one dimensional case, the associated stationary problem has only finite number of solutions (see \cite[Theorem 3.4]{Z86}). As a consequence, the connected set $\omega(\phi_0)$ consists only one point so that the global solution $\phi(t)$ to the time-dependent Cahn--Hilliard equation must converge to a certain equilibrium as $t \to +\infty$ (see \cite[Theorem 3.5]{Z86}). In this regard, we also refer to \cite{SZ93} for an extended result for the one dimensional non-isothermal Cahn--Hilliard system.

The situation in higher spatial dimensions is more complicated. Although the dissipative structure of problem \eqref{claCH}--\eqref{bdN2} (cf. \eqref{BELz}) guarantees that the $\omega$-limit set of a given initial datum $\phi_0$ contains only the steady states. However, this property is not sufficient for showing that $\omega(\phi_0)$ is a singleton. One difficulty is that, in higher dimensions, the study on stationary solutions to the Cahn--Hilliard equation turns out to be rather involved (see e.g., \cite{CK96,WW98a,WW98b}). In particular, due to the possible non-convexity of the free energy $E^{\mathrm{bulk}}(\phi)$, the structure of the set of equilibria seems far from being well understood and estimates on the number of equilibria are only available in the one dimensional case so far \cite{Z86,GN95}. On the other hand, examples have been given in the literature such that for certain simple semilinear second order parabolic equation with a specific smooth nonlinear term, it admits a globally bounded solution whose $\omega$-limit set is diffeomorphic to the unit circle $S^1$, namely, a continuum \cite{PP96,PS02}.

 Various assumptions have been made in the literature to assure that any bounded global solution of an evolution equation will converge to a single equilibrium as $t\to +\infty$. Among these attempts, an efficient approach was proposed by \cite{Si83}, in which the author generalized the \L ojasiewicz inequality for analytic functions in the finite-dimensional space $\mathbb{R}^m$ to the infinite dimensional spaces and proved that, if the nonlinear term is real analytic in the unknown function then the convergence to a single equilibrium holds.

For problem \eqref{claCH}--\eqref{bdN2} with a constant mobility ($M=1$) and a general regular potential $F:\Omega\times \mathbb{R}\to \mathbb{R}$ depending analytically on $x$, $\phi$, uniformly in $x$ with suitable growth assumptions, the authors of \cite{RH99} applied the {\L}ojasiewicz--Simon approach and proved that the unique global solution converges to a single equilibrium as $t\to +\infty$ in the topology of $C^{4,\mu}(\Omega)$, $\mu\in (0,1)$ (see \cite[Theorem 3.1]{RH99}). The proof relies on a generalized gradient inequality of {\L}ojasiewicz--Simon type (see \cite[Theorem 3.2]{RH99}), which is applicable to the $H^{-1}$-gradient flow defined by the Cahn--Hilliard equation (cf. \eqref{grad1}). Indeed, from the energy dissipation relation \eqref{BELz} (with $M=1$) and the equation \eqref{claCH}, we see that
$$
\int_0^{+\infty}\|\nabla \mu(t)\|^2_{L^2(\Omega)}\,\d t<+\infty\quad \Longrightarrow \int_0^{+\infty}\|\partial_t\phi(t)\|^2_{(H^1(\Omega))'}\,\d t<+\infty.
$$
On the other hand, from the precompactness of the trajectory $\phi(t)$ in certain space, for instance, $C^{4,\mu}(\Omega)$ as shown in \cite{RH99}, we can find a convergent subsequence $\phi(t_n)\to \phi_\infty$ in $C^{4,\mu}(\Omega)$, as $t_n\nearrow +\infty$ for some equilibrium $\phi_\infty\in \omega(\phi_0)$ (thanks to the characterization on $\omega(\phi_0)$ for gradient systems). At this stage, if we can show (which is not obvious, since $h(t)\in L^2(0,+\infty)$ does not imply $h(t)\in L^1(0,+\infty)$)
$$
\int_{t_0}^{+\infty}\|\partial_t\phi(t)\|_{(H^1(\Omega))'}\,\d t<+\infty,\quad \text{for some}\ t_0>0,
$$
then
$$
\lim_{t\to+\infty}\int_{t}^{+\infty}\|\partial_t\phi(\tau)\|_{(H^1(\Omega))'}\,\d \tau=0.
$$
As a consequence, it holds
\begin{align*}
\|\phi(t)-\phi_\infty\|_{(H^{1}(\Omega))'}
&\leq \|\phi(t)-\phi(t_n)\|_{(H^{1}(\Omega))'}+ \|\phi(t_n)-\phi_\infty\|_{(H^{1}(\Omega))'}\\
&\leq \int_{t_n}^t \|\partial_t\phi(\tau)\|_{(H^1(\Omega))'}\,\d \tau
+\|\phi(t_n)-\phi_\infty\|_{(H^{1}(\Omega))'}\\
&\to 0,\quad \text{for all}\ t\geq t_n,\ \text{as}\ t_n\to +\infty,
\end{align*}
which further implies the convergence of $\phi(t)$ to $\phi_\infty$ for all $t\to +\infty$, namely, $\omega(\phi_0)=\{\phi_\infty\}$. We emphasize that the {\L}ojasiewicz--Simon type inequality played an essential role in obtaining the required $L^1$-integrability of $\|\partial_t\phi(t)\|_{(H^1(\Omega))'}$ on $\mathbb{R}^+$ from its $L^2$-integrability on $\mathbb{R}^+$, where the latter is an easy consequence from the energy dissipation. To achieve such a goal (not limited to the Cahn--Hilliard equation but also for other gradient-like systems), different arguments have been developed in the literature, we may refer to, for instance, \cite{Chi,FS00,HP01,Si83,RH99,J98} and also the book \cite{Huang06}.

Later in \cite{AW07}, the authors studied problem \eqref{claCH}--\eqref{bdN2} with $M=1$ and a logarithmic type potential like \eqref{logpo}. They derived an extended \L ojasiewicz--Simon type inequality (see \cite[Proposition 6.3]{AW07}) and applied it to show that every global solution converges to a single equilibrium as $t\to+\infty$ in $H^{2r}(\Omega)$ for $r\in (0,1)$. In this case, we note that regularity of stationary solutions, in particular, the property of strict separation from pure states (which is indeed a direct consequence from the observation that every solution of the stationary Cahn--Hilliard equation solves the viscous Cahn--Hilliard equation), plays a crucial role in the derivation of the \L ojasiewicz--Simon inequality involving a singular potential. In a similar sprit, the authors of \cite{GG06} presented a direct constructive descent method for finding minimizers of nonconvex functionals via the \L ojasiewicz--Simon approach and applied it to phase separation problems in multicomponent systems as well as image segmentation. Recently, in \cite{GGM17} the authors extended the convergence result for global weak solutions to the Cahn--Hilliard--Oono equation \eqref{claCHO} with a logarithmic potential subject to \eqref{ini0}--\eqref{bdN2}, see \cite[Theorem 7.1]{GGM17}. We would like to mention that the above convergence results for problem \eqref{claCH}--\eqref{bdN2} can be extended to the case with non-constant mobilities, provided that $M(\phi)$ is non-degenerate. Furthermore, we refer to the recent works \cite{Ab09,GG10,GGW18,JWZ15,WW12,ZWH09} and the references cited therein for generalizations to some more complicated systems with fluid interaction.

\section{Cahn--Hilliard Equation with Dynamic Boundary Conditions}
\setcounter{equation}{0}

The influence of boundaries (solid walls) on the phase separation process of binary mixtures has attracted a lot of attentions of scientists. For instance, the occurring structures of a binary polymer mixture during the phase separation process may get frozen by a rapid quench into the glassy state and micro-structures at surfaces on small length scales can be produced \cite{FM97}. In order to describe the short-range interactions between the wall and both components of the mixture, suitable surface free energy functional should be introduced into the system (see e.g., \cite{BF91,FM97,FR98})
\begin{align}
E^\mathrm{surf}(\phi)=\int_{\partial \Omega} \frac{\kappa}{2} |\nabla_\Gamma \phi|^2+G(\phi) \, \d S,
\label{Esurf}
\end{align}
where $\nabla_\Gamma$ stands for the tangential (surface) gradient operator defined on the boundary $\partial\Omega$ and $G$ is a surface potential function that takes different forms in different physical contexts. The coefficient $\kappa\geq 0$ is related to the influence of spatial order parameter fluctuations on $\partial \Omega$ (e.g., the surface diffusion). When $\kappa=0$, the model is closely related to the evolution of a free interface in contact with the solid boundary, i.e., the well-known moving contact line problem \cite{Ja00,TR89}.

Like in Section 2, the Cahn--Hilliard equation for phase separation process with boundary effects still has to be supplemented by two boundary conditions, with natural considerations from the mass conservation and the energy dissipation (or in other words, energy balance). However, the difference is that one now has to take into account the non-trivial boundary effects driven by the surface energy \eqref{Esurf}. Different types of dynamic boundary conditions have been proposed and analyzed in the literature. In what follows, we shall review some recent progresses in this direction. Before proceeding, we just remark that dynamic boundary conditions can be found in various physical problems (e.g., heat conduction with heat source on the boundary, kinetic motion on the boundary of a vibrating object etc) and have been extensively analyzed in the literature, see for instance, \cite{Ai95,AP10,E92,EGG03,Fr10,FS68,Ga12,Go06,GS93,Hi89,La99,SW10,WH07a,XL04} and the references cited therein. In particular, we refer to \cite{CGGM10,CGM08,CM09,GG08,GG09,GMS10} for the Caginalp phase-field system with dynamic boundary conditions accounting for the non-isothermal phase transition process of a two-phase material.

\begin{remark}
In the remaining part of this section, since we do not consider asymptotic behavior with respect to the parameter $\epsilon$ in $E^\mathrm{bulk}$ (recall \eqref{Ebulk}), without loss of generality, we can simply set $\epsilon=1$.
\end{remark}

\subsection{Dynamic boundary condition of Allen--Cahn type}

 The first set of boundary conditions we would to mention are those proposed in \cite{FM97,K-etal} (see also \cite{BF91} for a derivation from the semi-infinite Ising model with Kawasaki spin exchange dynamics). To ensure the conservation of mass in the bulk (cf. \eqref{massz}), one can again assume the following no-flux boundary condition (for the sake of simplicity, hereafter we set the mobility $M=1$)
\begin{equation}
\partial_\mathbf{n} \mu =0, \quad\  \text{on}\ \partial\Omega \times (0,T).\label{bdN2da}
\end{equation}
On the other hand, to guarantee that the Cahn--Hilliard system tends to minimize its total free energy, the following variational boundary condition was proposed:
\begin{equation}
\frac{1}{\Gamma_s}\partial_t \phi -\kappa\Delta_\Gamma \phi +\partial_\mathbf{n}\phi +G'(\phi)=0, \quad  \text{on}\ \partial \Omega\times (0,T),
\label{bdN2d}
\end{equation}
where $\Gamma_s>0$ denotes surface kinetic coefficient, $\Delta_\Gamma$ stands for the Laplace--Beltrami operator on $\partial \Omega$ and
the interaction from the bulk is presented by the term $\partial_\mathbf{n}\phi$.
The above boundary condition is usually referred to as a dynamic one, since it contains the time derivative of the order parameter. It can be viewed as a relaxation dynamics ($L^2$-gradient flow) of the surface free energy $E^\mathrm{surf}$ on $\partial \Omega$, and moreover, it yields the energy dissipation relation
\begin{align}
&\frac{\d}{\d t} \Big[E^\mathrm{bulk}(\phi(t))+E^\mathrm{surf}(\phi(t))\Big]
+ \int_\Omega |\nabla \mu|^2 \,\d x
+ \frac{1}{\Gamma_s} \int_{\partial \Omega} |\partial_t \phi|^2\, \d S=0,\quad \forall\, t\in (0,T).
\label{bbela}
\end{align}
When $\kappa=0$, the dynamic boundary condition \eqref{bdN2d} was proposed to describe some other physically relevant situation for the dynamics of fluid–fluid free interface at the solid boundary (see e.g., \cite{Q1}) such that the contact angle on $\partial \Omega$ turns out to be time dependent and may deviate from the static one like $\pi/2$ (cf. \eqref{bdN2}).

Besides the Lyapunov structure \eqref{bbela}, the Cahn--Hilliard equation \eqref{claCH} subject to boundary conditions \eqref{bdN2da}--\eqref{bdN2d} can actually be interpreted as a gradient flow of the total energy $E^{\mathrm{total}}(\phi)=E^\mathrm{bulk}(\phi)+E^\mathrm{surf}(\phi)$ in the space $\dot{H}^{-1}(\Omega)\times L^2(\partial\Omega)$ with respect to the following inner product:
$$
(g,f)=\int_{\Omega} \nabla(\mathcal{N} g) \cdot \nabla (\mathcal{N} f) \, \d x + \int_{\partial\Omega} g f \,\d S.
$$
Roughly speaking, the gradient flow is of a mixed type such that it is $\dot{H}^{-1}$ in the bulk and $L^2$ on the boundary (see \cite[Section 3]{GK20}).\medskip

\textbf{Well-posedness and long-time behavior}. Now we review some presentive works on the well-posedness of the initial boundary value problem
\begin{equation}
\left\{
\begin{aligned}
&\partial_t \phi= \Delta \mu,\quad &&\text{in}\ \Omega \times (0,T),\\
&\mu=-\Delta \phi + F'(\phi),\quad &&\text{in}\ \Omega \times (0,T),\\
&\partial_\mathbf{n} \mu =0, \quad && \text{on}\ \partial\Omega \times (0,T),\\
&\frac{1}{\Gamma_s} \partial_t \phi-\kappa\Delta_\Gamma \phi +\partial_\mathbf{n}\phi +G'(\phi)=0, \quad && \text{on}\ \partial \Omega\times (0,T),\\
&\phi|_{t=0}=\phi_0(x),\quad &&  \text{in}\ \Omega.
\end{aligned}
\right.
\label{claCHdb1}
\end{equation}

Let us first focus on the case with $\kappa>0$. The first result on existence and uniqueness of global strong solutions to problem \eqref{claCHdb1} was obtained in \cite{RZ03}. There the authors considered the following specific choice of free energies
\begin{equation}
\begin{cases}
E^{\mathrm{bulk}}(\phi)=\displaystyle{\int_\Omega \left(\frac{1}{2}|\nabla \phi|^2+ \frac{1}{4}\phi^4-\frac12\phi^2\right)\d x,}\smallskip \\
E^\mathrm{surf}(\phi)=\displaystyle{\int_{\partial \Omega} \left(\frac{\kappa}{2} |\nabla_\Gamma \phi|^2+\frac{g_s}{2}\phi^2-h_s\phi\right) \, \d S,}
\end{cases}
\label{popo}
\end{equation}
where the constant $g_s > 0$ accounts for a modification of the effective interaction between the components at the wall, and $h_s\neq 0$ describes the possible preferential attraction of one of the two components by the wall \cite{K-etal}. To overcome the mathematical difficulties due to the presence of the dynamic boundary condition as well as the Laplace--Beltrami operator on the boundary, they introduce an approximate problem related to the phase-field system of Caginalp type
\begin{equation}
\left\{
\begin{aligned}
&\varepsilon \partial_t \mu -\Delta \mu=-\partial_t \phi,\quad &&\text{in}\ \Omega \times (0,T),\\
&\partial_\mathbf{n} \mu =0, \quad && \text{on}\ \partial\Omega \times (0,T),\\
&\varepsilon \partial_t \phi-\Delta \phi =\mu-F'(\phi),&&\text{in}\ \Omega \times (0,T),\\
&\frac{1}{\Gamma_s} \partial_t \phi-\kappa\Delta_\Gamma \phi +\partial_\mathbf{n}\phi +G'(\phi)=0, \quad && \text{on}\ \partial \Omega\times (0,T),\\
&\mu|_{t=0}=\mu_0(x),\ \ \phi|_{t=0}=\phi_0(x),\quad &&  \text{in}\ \Omega,
\end{aligned}
\right.
\label{claCHdb1C}
\end{equation}
for some $\varepsilon>0$.
Here, we note that \eqref{claCHdb1C} is slightly different from the exact form of the approximate problem in \cite[Section 2]{RZ03}. The reason is that we just want to mention the idea but not to get into some technical details therein. The approximate problem \eqref{claCHdb1C} can be locally solved by first analyzing its linearized problem and then applying the contraction mapping theorem \cite[Theorem 3.1]{RZ03}. After deriving uniform-in-time a priori estimates, they were able to extend the unique local solution to be a global one \cite[Theorem 4.1]{RZ03}. Finally, with the help of uniform a priori estimates that are independent of $\varepsilon$, they passed to the limit as $\varepsilon\to 0^+$, to obtain the existence of global strong solutions to the original problem \eqref{claCHdb1}.

Uniqueness can be proved by using the standard energy method, however, it was not clear whether the solution obtain in \cite{RZ03} defines a $C_0$-semigroup in the energy space (cf. \eqref{Ebulk}, \eqref{Esurf} and note that here $\kappa>0$):
\begin{equation}
\mathcal{V}^1=\big\{(\phi, \psi)\in H^1(\Omega)\times H^{1}(\partial \Omega):\ \psi=\phi|_{\partial\Omega} \big\}.
\label{V1}
\end{equation}
To handle this issue, in a subsequent paper \cite{PRZ06}, the authors studied the maximal $L^p$-regularity of the system \cite[Theorem 2.1]{PRZ06} and provided an alternative proof on the global well-posedness of problem \eqref{claCHdb1} (see \cite[Theorem 4.1]{PRZ06}), which further implies that the solution indeed defines a $C_0$-semigroup in $\mathcal{V}^1$ as expected. Based on this fact, they were able to prove the existence of a global attractor in suitable phase spaces \cite[Theorem 5.1, Theorem 5.2]{PRZ06}. Later, the maximal $L^p$-regularity approach was also used in \cite{PW06} to analyze the non-isothermal Cahn--Hilliard equation with dynamic boundary conditions. For further applications to parabolic problems with boundary dynamics of relaxation type, we refer to \cite{DPZ08}.

A third approach to study problem \eqref{claCHdb1} was given in \cite{MZ05}, where  the authors considered the viscous Cahn--Hilliard equation for $\alpha\geq 0$:
\begin{equation}
\left\{
\begin{aligned}
&\partial_t \phi= \Delta \mu,\quad &&\text{in}\ \Omega \times (0,T),\\
&\mu=\alpha\partial_t \phi -\Delta \phi + F'(\phi),\quad &&\text{in}\ \Omega \times (0,T),\\
&\partial_\mathbf{n} \mu =0, \quad && \text{on}\ \partial\Omega \times (0,T),\\
&\partial_t \phi-\kappa\Delta_\Gamma \phi +\partial_\mathbf{n}\phi +G'(\phi)=0, \quad && \text{on}\ \partial \Omega\times (0,T),\\
&\phi|_{t=0}=\phi_0(x),\quad &&  \text{in}\ \Omega,
\end{aligned}
\right.
\label{claCHdb1v}
\end{equation}
with some general (but regular) potentials $F,G\in C^3(\mathbb{R})$ with arbitrary growth and satisfying the following dissipative conditions:
$$
\liminf_{|s|\to+\infty}F''(s)>0, \quad \liminf_{|s|\to+\infty}G''(s)>0.
$$
They introduced a new variable on the boundary, i.e., by taking the trace of the phase function such that
$$
\psi=\phi|_{\partial\Omega}
$$
and then treated the dynamic boundary condition as a separate evolution equation (parabolic when $\kappa>0$) on
$\partial \Omega$, namely,
\begin{equation}
\left\{
\begin{aligned}
&\partial_t \phi= \Delta \mu,\quad &&\text{in}\ \Omega \times (0,T),\\
&\mu=\alpha\partial_t \phi -\Delta \phi + F'(\phi),\quad &&\text{in}\ \Omega \times (0,T),\\
&\partial_\mathbf{n} \mu =0, \quad && \text{on}\ \partial\Omega \times (0,T),\\
&\phi=\psi,&& \text{on}\ \partial\Omega \times (0,T),\\
&\partial_t \psi-\kappa\Delta_\Gamma \psi +\partial_\mathbf{n}\phi +G'(\psi)=0, \quad && \text{on}\ \partial \Omega\times (0,T),\\
&\phi|_{t=0}=\phi_0(x),\quad &&  \text{in}\ \Omega,\\
&\psi|_{t=0}=\psi_0(x):=\phi_0(x)|_{\partial\Omega},\quad && \text{on}\ \partial \Omega.
\end{aligned}
\right.
\label{claCHdb1v2}
\end{equation}
Within this convenient framework, the authors proved the existence and uniqueness of global solutions to problem \eqref{claCHdb1v2} by the Leray--Schauder principle, see  \cite[Theorem 2.1, Corollary 2.2, Theorem 2.2]{MZ05}. Moreover, they constructed a robust (as $\alpha\to 0^+$) family of exponential attractors for the semigroups associated with problem \eqref{claCHdb1v2} in suitable phase spaces with constraints due to the mass conservation \cite[Theorem 3.1]{MZ05}. As a byproduct, the result therein also implies that the global attractor obtained in \cite{PRZ06} has finite fractal dimension.

The analysis for problem \eqref{claCHdb1} with singular potentials are more involved, since the combination of dynamic boundary conditions and of singular potentials can produce additional strong singularities on the corresponding solutions close to the boundary \cite{MZ10}. In this direction, the first result on existence and uniqueness of global weak solutions was obtained in \cite{GMS09}, within the following general setting
\begin{equation}
\left\{
\begin{aligned}
&\partial_t \phi= \Delta \mu,\quad &&\text{in}\ \Omega \times (0,T),\\
&\mu=\alpha\partial_t \phi -\Delta \phi + \beta(\phi)+\pi(\phi)-f,\quad &&\text{in}\ \Omega \times (0,T),\\
&\partial_\mathbf{n} \mu =0, \quad && \text{on}\ \partial\Omega \times (0,T),\\
&\phi=\psi,&& \text{on}\ \partial\Omega \times (0,T),\\
&\partial_t \psi-\kappa\Delta_\Gamma \psi +\partial_\mathbf{n}\phi +\beta_\Gamma(\psi)+\pi_\Gamma(\psi)=f_\Gamma, \quad && \text{on}\ \partial \Omega\times (0,T),\\
&\phi|_{t=0}=\phi_0(x),\quad &&  \text{in}\ \Omega,\\
&\psi|_{t=0}=\psi_0(x):=\phi_0(x)|_{\partial\Omega},\quad && \text{on}\ \partial \Omega.
\end{aligned}
\right.
\label{claCHdb1v2g}
\end{equation}
where $\alpha, \kappa\geq 0$, see \cite[Theorem 1, Theorem 2]{GMS09}.
In the above system, the authors of \cite{GMS09} made decompositions of the functions $F'$ and $G'$ as $F' = \beta + \pi$ and $G' = \beta_\Gamma + \pi_\Gamma$, respectively, where $\beta$ and $\beta_\Gamma$ are monotone and possibly non-smooth, while $\pi$ and $\pi_\Gamma$ are some regular (Lipschitz) perturbations. In \cite{GMS09}, they aimed to keep the form of nonlinearities as general as possible. For instance, $\beta$ and $\beta_\Gamma$ are allowed to be essentially arbitrary, with some compatibility conditions such that $\beta$ grows faster than $\beta_\Gamma$ (i.e., the bulk potential plays a dominating role) and that the other boundary contributions satisfy a specific sign condition (see \cite[Remark 6]{GMS09}). Well-posedness results under alternative assumptions for $\beta$, $\beta_\Gamma$ like some growth conditions were also discussed in \cite[Theorem 3]{GMS09}. The proofs for the existence results in \cite{GMS09} rely on the Yosida regularization of the possibly singular potentials, a suitable Faedo--Galerkin scheme (via eigenfunctions of the elliptic problem $-\Delta u = \lambda u$ subject to homogeneous Neumann boundary condition) together with the compactness method (based on suitable a priori estimates performed on the approximate solutions). This argument indeed provides the fourth approach to handle problem \eqref{claCHdb1} in a general setting.

Some further investigation was performed in \cite{MZ10}. We recall that the result in \cite{GMS09} yields the existence of global weak solutions if the (regular) surface nonlinearity $G'$ has the ``right" sign at the singular points (i.e., pure states) of the bulk nonlinearity $F'$, that is, $\pm G'(\pm 1)>0$. However, when the sign condition is violated, it was shown in \cite[Remark 6.2]{MZ10} that one may not have classical solutions even in the one  dimensional case. To overcome this difficulty,  the authors of \cite{MZ10} introduced a suitable notion of variational solutions (see \cite[Definition 3.1]{MZ10}) and proved its existence and uniqueness by using proper approximations of the singular potentials (see \cite[Theorem 3.2]{MZ10}). Then they proved the existence of global and exponential attractors in \cite[Theorem 5.2]{MZ10}. Besides, connections between the variational solution and the usual distribution solution were established in \cite{MZ10}. The authors discussed the possible separation of the solutions from the singularities of the singular bulk potential (see \cite[Section 4]{MZ10}). They showed that those variational solutions are H\"{o}lder continuous in space and will become solutions in the usual sense if they do not reach the pure states on the boundary, which can be guaranteed if either the sign condition as in \cite{GMS09} holds or the singularity of the bulk potential $F$ is strong enough (unfortunately, not satisfied by the logarithmic one \eqref{logpo}). We refer to \cite{CMZ11,Mi17,Mi19,MZ10} for detailed discussions.

In \cite{GMS09,MZ10} the authors mainly treated the case that the bulk potential $F$ plays a dominating role.
Motivated by the study for the Allen--Cahn equation with singular potentials and dynamic boundary condition in \cite{CC13}, the authors of \cite{CGS14} considered the opposite side of the compatibility condition for potential functions such that the boundary potential $G$ is now the leading one (see \cite[(2.11)]{CGS14}). For problem \eqref{claCHdb1v} with $\alpha \geq 0$ and some general assumptions on the potentials that cover all the typical cases \eqref{logpo}, \eqref{regpo} and \eqref{obspo}, the authors proved global existence, uniqueness and regularity results on its solutions (see \cite[Theorem 2.2, Theorem 2.3, Theorem 2.4, Theorem 2.6]{CGS14}). The argument therein implies that in the case of a dominating boundary potential, the analysis can be somewhat simplified and it indeed allows for a unified treatment of the initial boundary value problem \eqref{claCHdb1v}. As far as some extensions are concerned, we refer to \cite{CF15} for the Cahn--Hilliard system with dynamic boundary condition and mass constraint on the boundary, to \cite{CGS17} for the study of a nonstandard viscous Cahn--Hilliard system with dynamic boundary condition, and to \cite{Sc19} for the well-posedness of a Cahn--Hilliard system with nonlinear viscosity terms and nonlinear dynamic boundary conditions.

When $\kappa=0$, i.e., the surface diffusion on $\partial \Omega$ is absent, we recall that existence and uniqueness of global weak solutions to problem \eqref{claCHdb1} have already been obtained in \cite{GMS09}. We also refer to \cite{CWX14} for analysis of problem \eqref{claCHdb1} in a different context, which accounts for the evolution of fluid-fluid free interfaces along the solid boundary in the ``slow" dynamics such that the effect of flow can be neglected. The authors proved well-posedness of the system with regular potentials \cite[Theorem 2, Theorem 3]{CWX14} and investigated its sharp interface limit via the method of matched asymptotic expansion. Besides, the dynamics of the contact point and the contact angle were described and the results were compared with numerical simulations.

Then a natural question arises, when the surface diffusion in the dynamic boundary condition is vanishing as $\kappa\to 0^+$, whether the solutions to the problem \eqref{claCHdb1} with surface diffusion will converge to the solutions of problem \eqref{claCHdb1} without surface diffusion? In \cite{CF20b}, the authors gave an affirmative answer for a wide class of potentials including \eqref{logpo}, \eqref{regpo} and \eqref{obspo} (see \cite[Theorem 2.2]{CF20b}). Their result indicates that the solution of the limiting problem with $\kappa=0$ will lose some spatial regularity due to the absence of the surface diffusion on $\partial\Omega$ as expected. For instance, the normal derivative  $\partial_\mathbf{n}\phi$ on the boundary in general is no longer a function but an element in a dual space. Only under some specific conditions it (as well as the boundary condition \eqref{bdN2d}) can hold almost everywhere (see \cite[Theorem 5.1]{CF20b}).

When the long-time behavior of global solutions to problem \eqref{claCHdb1} is concerned, we refer to \cite{GMS09a,MZ05,MZ10,Mi17,Mi19,PRZ06} for extensive studies within the theory of global and exponential attractors, see also \cite{CGG,CGW,Ga07,Ga08,GG13,GM1,GM2} for related results on some extended systems. The role of surface diffusion in dynamic boundary conditions for parabolic and elliptic equations was also discussed in \cite{Ga15} within the attractor theory.

Thanks to the energy dissipation relation \eqref{bbela}, we see that the dynamics of the system can be viewed as a mixing of the Cahn--Hilliard type in the bulk combined with the Allen--Cahn type on the boundary. In particular, we expect that the global solution of problem \eqref{claCHdb1} will converge to a single equilibrium as $t\to+\infty$ under suitable assumptions on the potentials $F$ and $G$. When $F$, $G$ take the specific form in \eqref{popo}, in \cite{WZ04}, the authors proved the convergence result (see \cite[Theorem 1.1]{WZ04}) by applying an extended \L ojasiewicz--Simon inequality that involves a boundary term \cite[Lemma 3.4]{WZ04}. In \cite{CFP06}, the authors obtained the convergence of global solutions to problem \eqref{claCHdb1} with a general regular potential $F$ being analytic and under suitable growth assumptions (see \cite[Theorem 2.3]{CFP06}). They derived a \L ojasiewicz--Simon type inequality in the dual space \cite[Proposition 6.6]{CFP06} by adapting the abstract result in \cite[Corollary 3.11]{Chi}. Moreover, they proved the convergence to a single equilibrium of global solutions to the phase-field model of Caginalp type with dynamic boundary conditions (cf. \eqref{claCHdb1C}), see also \cite{WGZ07} for the convergence result of an extended phase-field system with hyperbolic relaxation but without boundary diffusion. More general case was treated in \cite{GMS09a}, where the authors considered a general singular bulk potential (including \eqref{logpo}) together with a regular boundary potential, in the framework of \cite{GMS09}. They extended the argument in \cite{CFP06} and derived an extended \L ojasiewicz--Simon inequality that also allows nonlinear term on the boundary \cite[Proposition 6.1]{GMS09a} (recall that in \cite{WZ04,CFP06} only linear boundary was considered). Then they proved the convergence to a single equilibrium for global solutions as $t\to+\infty$.
For related results on some extended Cahn--Hilliard systems, we may refer to \cite{CGG,CGW,GM2,GGM16} and the references cited therein.

We also mention \cite{Ga17} for the study on a nonlocal version of the Cahn--Hilliard equation characterized by a fractional diffusion operator which is subject to a fractional dynamic boundary conditions. For the case with regular potentials, the author proved global well-posedness, regularity of solutions \cite[Theorems 3.3, 3.4, 3.5]{Ga17} and the existence of an exponential attractor \cite[Theorem 4.3]{Ga17}, which yields existence of the global attractor with finite fractal dimension \cite[Corollary 1]{Ga17}.

\subsection{Dynamic boundary condition of Cahn--Hilliard type (I)}
For simplicity, hereafter we consider the Cahn--Hilliard equation \eqref{claCH} with $M=1$. In \cite{Ga06}, the author considered the phase separation in a binary mixture confined to a bounded domain with porous walls, which may be (semi) permeable. In this case, due to the possible mass transfer on and through the boundary $\partial\Omega$, the homogeneous Neumann boundary condition \eqref{bdN2da} for the chemical potential should be modified.
Alternatively, the author of \cite{Ga06} proposed the following Wentzell type boundary condition for $\mu$
\begin{equation}
\Delta \mu + b \partial_\mathbf{n}\mu+ c\mu=0,\quad \text{on}\ \partial\Omega \times(0,T),
\label{bdWe}
\end{equation}
for some $b>0, c\geq 0$ (cf. \cite{FGGR02,Go06} for the Wentzell boundary condition for the heat equation and wave equation). Then for sufficiently smooth solutions $(\phi, \mu)$, \eqref{bdWe} is equivalent to the following dynamic boundary condition
\begin{equation}
\partial_t \phi + b \partial_\mathbf{n}\mu+ c\mu=0,\quad \text{on}\ \partial\Omega \times(0,T).
\label{bdWe1}
\end{equation}
Under the above choice, we no longer have mass conservation in the bulk (unless letting $b\to+\infty$), but a new conservation law for the ``total mass" defined in the measure space $(\overline{\Omega},\d\nu)=(\Omega,\d x)\oplus(\partial\Omega,\d S/b)$ (see \cite{Ga06,Go06}):
\begin{equation}
\frac{\d}{\d t}\left(\int_\Omega\phi\,\d x+ \int_{\partial\Omega} \phi\, \frac{\d S}{b}\right)= -\int_{\partial\Omega} c\mu\,\frac{\d S}{b}.\nonumber
\end{equation}
The boundary condition \eqref{bdWe} describes the situation that there is certain
mass source (or leak) on the boundary and thus the system undergoes a different diffusive process. When $c=0$, we see that the total mass is conserved along the evolution such that
\begin{equation}
\int_\Omega \phi(t) \,\d x+ \int_{\partial\Omega} \phi(t) \,\frac{\d S}{b} =\int_\Omega \phi(0)\,\d x+\int_{\partial\Omega} \phi(0) \,\frac{\d S}{b}, \quad \forall\, t\in [0,T].
\label{massW}
\end{equation}
Namely, the wall is non-permeable, but certain mass exchange between the bulk and the boundary is allowed. On the other hand, in order to guarantee that the system tends to minimize its total free energy, the author proposed the following variational boundary condition in \cite{Ga06} (with $G'$ being a linear function of $\phi$):
\begin{equation}
-\kappa \Delta_\Gamma\phi+\partial_\mathbf{n}\phi + G'(\phi) =\frac{\mu}{b},\quad \text{on}\ \partial\Omega \times(0,T).
\label{bdWe2}
\end{equation}
As a consequence, it holds
\begin{equation}
\frac{\d}{\d t} \Big[E^\mathrm{bulk}(\phi(t))+E^\mathrm{surf}(\phi(t))\Big]  + \int_\Omega |\nabla \mu|^2\,\d x + c \int_{\partial \Omega} |\mu|^2\, \frac{\d S}{b}=0,\quad \forall\, t\in (0,T).
\label{bbelbW}
\end{equation}

Later in \cite{GMS11}, the authors considered the phase separation process in a bounded domain with non-permeable walls in some more general setting. Requiring that the total free energy is decreasing in time and the total mass is conserved (see \eqref{massW}), they derived the following set of boundary conditions via a variational principle (see also \cite{CGM12}),
\begin{align}
&\partial_t\phi +b( \partial_\mathbf{n} \mu - \sigma \Delta_\Gamma \mu) =0, && \text{on}\ \partial\Omega \times (0,T),\label{bdN1a}\\
&\frac{\mu}{b}=-\kappa \Delta_\Gamma \phi +\partial_\mathbf{n}\phi +G'(\phi),&& \text{on}\ \partial\Omega \times (0,T),\label{bdN2a}
\end{align}
for some $b>0$, $\sigma\geq 0$, $\kappa\geq 0$. Here, the term $\Delta_\Gamma \mu$ in \eqref{bdN1a} corresponds to some regularizing effect for $\mu$ on the boundary. For $\sigma>0$ and $\kappa>0$, \eqref{bdN1a} is often referred to as the Cahn--Hilliard type dynamic boundary condition in the literature. When $\sigma=0$, \eqref{bdN1a} simply reduces to the Wentzell boundary condition \eqref{bdWe1} with $c=0$. In \cite{GMS11}, the authors indeed considered a more general situation such that the parameter $b$ is replaced by a positive function $w^{-1}$: $w\in L^\infty(\partial\Omega)$, $0<w_*\leq w(x)\leq w^*$ for a.e. $x\in \partial\Omega$, where $w_*, w^*$ are given constants. With this generalization, they were able to handle the situation when the ``boundary mass" is linked to the variable $\phi$ in a different way with respect to the ``bulk mass", for instance, the case when the dynamic boundary condition arises as an approximation of a thin diffusive layer occupied by a different material (see \cite[Section 1]{GMS11}). From \eqref{bdN1a} and \eqref{bdN2a}, it follows that the mass conservation property  \eqref{massW} is satisfied and the following energy dissipation relation holds (cf. \eqref{bbelbW})
\begin{align}
&\frac{\d}{\d t} \Big[E^\mathrm{bulk}(\phi(t))+E^\mathrm{surf}(\phi(t))\Big]  + \int_\Omega |\nabla \mu|^2\,\d x + \sigma \int_{\partial \Omega} |\nabla_\Gamma \mu|^2\, \frac{\d S}{b}=0,\quad \forall\, t\in (0,T).
\label{bbelb}
\end{align}

\textbf{Gradient flow structure}.
As pointed out in \cite{GK20}, the Cahn--Hilliard equation \eqref{claCH} (with $M=1$) subject to \eqref{bdN1a}--\eqref{bdN2a} admits a gradient flow structure. For simplicity, set $b=1, \sigma>0$. Let us consider the elliptic problem of Wentzell Laplacian
\begin{equation*}
\begin{cases}
-\Delta u=f_\Omega,\qquad\qquad\quad  \text{in}\ \Omega,\\
-\sigma\Delta_\Gamma u +\partial_\mathbf{n} u=f_\Gamma, \quad\, \text{on}\ \partial\Omega.
\end{cases}
\end{equation*}
Applying the Lax--Milgram theorem, one can see that the above problem admits a unique weak solution $u = \mathcal{W} (f)\in \mathcal{V}^1$ with constraint $\int_\Omega u\,\d x+\int_{\partial\Omega} u\,\d S = 0$, if the two terms on the right-hand
side satisfy $f=(f_\Omega,f_\Gamma)\in (\mathcal{V}^1)'$ and $ \langle f_\Omega, 1\rangle_{(H^1(\Omega))',H^1(\Omega)}+\langle f_\Gamma, 1\rangle_{(H^1(\partial\Omega))',H^1(\partial\Omega)}= 0$. With the help of the solution operator $\mathcal{W}$, equation \eqref{claCH} subject to \eqref{bdN1a}--\eqref{bdN2a} can be viewed as a gradient flow of the total energy $E^{\text{total}}(\phi)$ with respect to the inner product
$$
(g,f)_{(\mathcal{V}^1)'}=\int_{\Omega} \nabla(\mathcal{W} g) \cdot \nabla (\mathcal{W} f) \, \d x + \sigma \int_{\partial\Omega} \nabla_\Gamma (\mathcal{W} g)\cdot \nabla_\Gamma(\mathcal{W} f) \,\d S.
$$
Similar conclusions could be drawn for general choices of parameters. \medskip

\textbf{Well-posedness and long-time behavior}. Assume that $\Omega\subset \mathbb{R}^d$, $d=1,2,3$, is a smooth bounded domain. Without loss of generality, let us write the resulting initial boundary value problem in the following form:
\begin{equation}
\left\{
\begin{aligned}
&\partial_t \phi= \Delta \mu,\quad &&\text{in}\ \Omega \times (0,T),\\
&\mu=-\Delta \phi + F'(\phi),\quad &&\text{in}\ \Omega \times (0,T),\\
&\partial_t \phi + b\partial_\mathbf{n} \mu - \sigma \Delta_\Gamma \mu + c \mu =0, \quad && \text{on}\ \partial\Omega \times (0,T),\\
&\frac{\mu}{b}=-\kappa\Delta_\Gamma \phi +\partial_\mathbf{n}\phi +G'(\phi), \quad && \text{on}\ \partial \Omega\times (0,T),\\
&\phi|_{t=0}=\phi_0(x),\quad &&  \text{in}\ \Omega.
\end{aligned}
\right.
\label{claCHdb1W}
\end{equation}
Note that in the one dimensional case, the terms involving $\Delta_\Gamma$ do not appear in any of the boundary conditions.

We first review some known results for the case $\sigma=0$.  The first results on existence and uniqueness of global strong solutions to problem \eqref{claCHdb1W} with a regular bulk potential $F$ and a quadratic boundary potential $G$ was obtained in \cite{Ga06} by employing some classical methods such as the contraction mapping principle and energy methods (see \cite[Theorem 11, Theorem 12]{Ga06}). The proof therein extends the idea in \cite{RZ03} by investigating an approximate problem of Caginalp type. Besides, a comparison between two solutions of the systems \eqref{claCHdb1} and \eqref{claCHdb1W} was made in  \cite[Theorem 14]{Ga06}, which provides an estimate on the difference between the two solutions on finite time intervals, explicitly in terms of the parameters $b$ and $\Gamma_s$. Later in \cite{Ga06a}, the author improved his result by using a different approach, i.e., the Faedo--Galerkin method, to prove the existence and uniqueness of a global solution to the same problem, but under some more general assumptions on $F$, see \cite[Theorem 4.1]{Ga06a}. The result implies that the solution defines a continuous semigroup on a suitable phase space, which enables him to investigate its long-time behavior by proving the existence of an exponential attractor (and thus of a global attractor with finite dimension), see \cite[Theorem 5.1, Theorem 5.6]{Ga06a}. In \cite{Ga08a}, the author further studied the asymptotic behavior as $b\to+\infty$ and constructed a robust family of exponential attractors for the problem under the same assumptions on the potential as in \cite{Ga06a}.

Convergence of global solutions to a single equilibrium as $t \to +\infty$ for problem \eqref{claCHdb1W} with $\sigma=0$ was first analyzed in \cite{WH07}. Under the assumption that $c, \kappa>0$ and $F$ is analytic with respect to the phase function, the author proved the expected convergence result \cite[Theorem 1.1]{WH07} by means of an extended \L ojasiewicz--Simon type inequality with boundary term (see \cite[Lemma 3.5]{WH07}). Higher order estimates on the convergence rate were also obtained by the \L ojasiewicz--Simon approach (cf. \cite{HJ01}) together with the energy method. The conserved case $c=0$ is somewhat more involved, since the boundary contribution in the energy dissipation of \eqref{bbelbW} vanishes. To overcome the difficulty coming from a weaker dissipation, in \cite{GaW08} the authors made use of a generalized Poincar\'{e} inequality to derive an improved \L ojasiewicz--Simon type inequality subject to certain mass constraint (see \cite[Lemma 4.1]{GaW08}), and then deduced the convergence result \cite[Theorem 2.4]{GaW08}. Moreover, they could prove the convergence to steady states of global solutions in the limiting case $b=\infty$.

Next, the case $\sigma>0$, $c=0$ was studied in \cite{GMS11} in a rather general setting, i.e., the bulk potential $F$ can be a singular one in a wide class and the constant $b$ can be a bounded positive function. There the authors overcome difficulties similar to those encountered in \cite{MZ10} and proved existence, uniqueness as well as regularity of global weak solutions (defined in a suitable sense) and studied their long-time behavior, including the existence of a compact global attractor and convergence to a single equilibrium as $t\to +\infty$.
The arguments in \cite{GMS11} were rather involved, since the authors tried to deal with the most general conditions on the nonlinear terms, see \cite[Section 3]{GMS11} for detailed discussions. Concerning a different approach, we refer to \cite{NK17} for existence and uniqueness of solutions to problem \eqref{claCHdb1W} with regular potentials on a domain with either permeable or non-permeable walls, which were obtained by applying the general theory of maximal $L_p$ regularity of relaxation type in \cite{DPZ08}.

We note that the boundary nonlinearity $G$ in \cite{GMS11} was assumed to be a regular function. The case with a possibly singular and dominating boundary potential was considered in \cite{CF15a}. The authors introduced the unknowns on the boundary
$$
\psi=\phi|_{\partial \Omega},\quad \mu_\Gamma=\mu|_{\partial\Omega}
$$
and considered the following alternative formulation of \eqref{claCHdb1W} (with some specific choice of parameters)
\begin{equation}
\left\{
\begin{aligned}
&\partial_t \phi= \Delta \mu,\quad &&\text{in}\ \Omega \times (0,T),\\
&\mu= - \Delta \phi + F'(\phi)-f,\quad &&\text{in}\ \Omega \times (0,T),\\
&\phi=\psi,\quad \mu=\mu_\Gamma && \text{on}\ \partial\Omega \times (0,T),\\
&\partial_t \psi +\partial_\mathbf{n} \mu -  \Delta_\Gamma \mu_\Gamma =0, \quad && \text{on}\ \partial \Omega\times (0,T),\\
&\mu_\Gamma= - \Delta_\Gamma \psi +\partial_\mathbf{n}\phi +G'(\psi)-f_\Gamma, \quad && \text{on}\ \partial \Omega\times (0,T),\\
&\phi|_{t=0}=\phi_0(x),\quad &&  \text{in}\ \Omega,\\
&\psi|_{t=0}=\psi_0(x):=\phi_0(x)|_{\partial\Omega},\quad && \text{on}\ \partial \Omega,
\end{aligned}
\right.
\label{claCHdb1Wg}
\end{equation}
where $f, f_\Gamma$ are some given functions. Working with a general setting of bulk and boundary potentials including the typical types \eqref{logpo}, \eqref{regpo} and \eqref{obspo}, they proved existence and uniqueness (indicated by a continuous dependence estimate) of global weak solutions \cite[Theorem 2.1, Theorem 2.2]{CF15a}. Refined regularity for weak solutions and existence of strong solutions were obtained as well, see \cite[Theorem 4.2]{CF15a}.
The proof of the existence result relies on the Yosida regularization of maximal monotone graphs and the introduction of viscous terms in the equations for both $\mu$, $\mu_\Gamma$ in \eqref{claCHdb1Wg}. Then the solvability of the approximate problem can be deduced from the abstract theory of doubly nonlinear evolution inclusions \cite{CV90}. After obtaining uniform estimates for the approximate solutions, the authors of \cite{CF15a} were able to pass to the limit and conclude the existence of weak solutions. See \cite{CF20a} for an extension to the Cahn--Hilliard equation on the boundary with bulk condition of Allen--Cahn type. We also refer to \cite{CGS18,CGS18a} for extended results on well-posedness and long-time behavior for the viscous regularization of problem \eqref{claCHdb1W} with an additional convection term. The proof for the existence of solutions in \cite{CGS18} is based on regularization of maximal monotone graphs but a different way of approximation for the equation via the Faedo--Galerkin method.

In the recent work \cite{FW21}, for problem \eqref{claCHdb1W} with $\sigma>0$, $c=0$ and
physically relevant singular (e.g., logarithmic) potentials, global regularity of weak solutions was established, see \cite[Theorem 2.1]{FW21}. When the spatial dimension is two, the authors showed the instantaneous strict separation property such that for arbitrary positive time any weak solution stays away from the pure phases $\pm 1$, while in the three dimensional case, although the instantaneous separation property remains open, an eventual separation property for large time was derived. As a consequence, they proved that every global weak solution converges to a single equilibrium as $t\to+\infty$ (see \cite[Theorem 2.2]{FW21}), by applying a suitable \L ojasiewicz--Simon type inequality. We expect that a similar result could be obtained for the case $\sigma=0$, $c\geq 0$.

\subsection{Dynamic boundary condition of Cahn--Hilliard type (II)}
From the previous discussions, we see that under all choices of boundary conditions (\eqref{bdN2da} with \eqref{bdN2d}, \eqref{bdWe1} with \eqref{bdWe2}, or \eqref{bdN1a} with \eqref{bdN2a}), the Cahn--Hilliard equation \eqref{claCH} satisfies two basic physical constraints, that is, the mass conservation and energy dissipation. Among them, the boundary conditions \eqref{bdN2da}, \eqref{bdWe1} and  \eqref{bdN1a} are proposed to keep suitable mass conservation property in the physical domain (see \eqref{massz}, \eqref{massW}), while the so-called variational boundary conditions \eqref{bdN2d}, \eqref{bdWe2}, \eqref{bdN2a} are chosen in a phenomenological way so that the validity of some specific energy dissipation relation is guaranteed (see \eqref{bbela}, \eqref{bbelbW}, or \eqref{bbelb}). We note that \eqref{bdN2d}, \eqref{bdWe2} and \eqref{bdN2a} can be viewed as some sufficient conditions for energy dissipation of the system, however, such choice may not be unique.

In \cite{LW}, the authors derived a different type of dynamic boundary condition for the Cahn--Hilliard equation \eqref{claCH}, that is,
\begin{equation}
\left\{
\begin{aligned}
&\partial_\mathbf{n} \mu=0, \quad &\text{on}\ \partial \Omega \times (0,T),\\
&\phi_t=\Delta_\Gamma (- \kappa \Delta_\Gamma \phi +\partial_\mathbf{n}\phi + G'(\phi)),\quad
&\text{on}\ \partial\Omega \times (0,T).
\end{aligned}
\right.
\label{bdCH}
\end{equation}
The derivation is based on an energetic variational approach that combines the least action principle and Onsager's principle of maximum energy dissipation, from which we see that the equation \eqref{claCH} subject to \eqref{bdCH} naturally fulfills three basic physical properties: (1) kinematics: conservation of mass both in the bulk $\Omega$ and on the boundary $\partial \Omega$; (2) energetics: dissipation of the total free energy; (3) force balance: both in the bulk $\Omega$ and on the boundary $\partial\Omega$.

Below we sketch the derivation of \eqref{bdCH}, and refer to \cite{LW} for further details in a more general setting. In the bulk $\Omega$, $\phi$ is assumed to be a locally conserved quantity that satisfies the continuity equation
\begin{equation}
\phi_t+\nabla\cdot(\phi \mathbf{u})=0,
\quad (x, t)\in \Omega\times (0,T),
\label{con1}
\end{equation}
where $\mathbf{u}: \Omega \to \mathbb{R}^d$ stands for the microscopic effective velocity (e.g., due to diffusion process etc).
We assume that $\mathbf{u}$  satisfies the no-flux boundary condition
\begin{equation}
\mathbf{u}\cdot \mathbf{n}=0,
\quad (x, t)\in \partial\Omega \times (0,T).
\label{bd1}
\end{equation}
The boundary dynamics is assumed to satisfy a local mass conservation law analogous to \eqref{con1} such that
\begin{equation}
\phi_t+\nabla_{\Gamma}\cdot (\phi \mathbf{v})=0,
\quad (x, t)\in \partial\Omega \times (0,T),
\label{con2}
\end{equation}
where $\mathbf{v}: \partial\Omega  \to \mathbb{R}^d$ denotes the microscopic effective tangential velocity field on the boundary.
We note that there is no need to impose any boundary condition on $\mathbf{v}$, since here the boundary $\partial \Omega$ is assumed to be a closed manifold.

Next, for an isothermal closed system, evolution of the binary mixtures is assumed to satisfy the following energy dissipation law
\begin{align}
 \frac{\d}{\d t}E^{\mathrm{total}}(t)=-\mathcal{D}^{\mathrm{total}}(t),\quad t\in(0,T),
 \label{ABEL}
\end{align}
where
\begin{equation}
E^{\mathrm{total}}(t)= E^{\mathrm{bulk}}(t)+E^{\mathrm{surf}}(t).
\label{fenergy}
\end{equation}
For instance, $E^{\mathrm{bulk}}$ and $E^{\mathrm{surf}}$ are given by \eqref{Ebulk} and \eqref{Esurf}, respectively.
On the other hand, the rate of energy dissipation $\mathcal{D}^{\mathrm{total}}$ is chosen as
\begin{align}
\mathcal{D}^{\mathrm{total}}(t)
=\mathcal{D}^{\mathrm{bulk}}(t)+\mathcal{D}^{\mathrm{surf}}(t),
\label{diss}
\end{align}
which also consists of contributions from the bulk and the boundary. Here, we assume that
\begin{align}
\mathcal{D}^{\mathrm{bulk}}(t)= \int_\Omega M_\mathrm{b}^{-1}\phi^2 \mathbf{u}\cdot \mathbf{u} \, \d x,
\quad \mathcal{D}^{\mathrm{surf}}(t)= \int_{\partial \Omega} M_\mathrm{s}^{-1} \phi^2  \mathbf{v}\cdot \mathbf{v} \, \d S,
\label{diss1}
\end{align}
where $M_\mathrm{b}$, $M_\mathrm{s}$ are some positive mobility functions.

Finally, in order to derive a closed system of partial differential equations, it remains to determine the microscopic velocities $\mathbf{u}$, $\mathbf{v}$ in equations \eqref{con1} and \eqref{con2}.
Using the principle of least action (LAP) and Onsager's maximum dissipation principle (MDP), we are able to derive the conservative and dissipative forces according to the free energy \eqref{fenergy} and the dissipation functional \eqref{diss}. Then by the force balance relation (Newton's second law), we obtain
\begin{equation}
\left\{
\begin{aligned}
&\phi \nabla \mu +M_\mathrm{b}^{-1}\phi^2 \mathbf{u}=0,\quad  &&\text{in}\ \Omega\times(0,T),\\
&\phi \nabla_\Gamma \left(\mu_\Gamma+ \partial_\mathbf{n} \phi\right)
+M_\mathrm{s}^{-1} \phi^2\mathbf{v}=0,\quad  &&\text{on}\ \partial\Omega\times(0,T),
\end{aligned}
\right.
\label{fba}
\end{equation}
where the chemical potentials are given by
\begin{align}
& \mu=-\Delta \phi +F'(\phi),\quad \mu_\Gamma=-\kappa \Delta_\Gamma \phi +\partial_\mathbf{n}\phi + G'(\phi).\nonumber
\end{align}
Solving $\mathbf{u}$, $\mathbf{v}$ from \eqref{fba} and inserting them back into \eqref{con1}, \eqref{con2}, we arrive at the following Cahn--Hilliard system subject to a new class of dynamic boundary condition:
\begin{equation}
\left\{
\begin{aligned}
&\partial_t \phi=\nabla \cdot(M_{\text{b}}\nabla \mu),\quad
&&\text{in}\  \Omega\times (0,T),\\
& \mu= - \Delta \phi +F'(\phi),\quad
&&\text{in}\  \Omega\times (0,T),\\
&\partial_\mathbf{n} \mu=0,\quad
&&\text{on}\ \partial \Omega \times (0,T),\\
&\partial_t \phi=\nabla_\Gamma\cdot (M_{\text{s}}\nabla_\Gamma \mu_\Gamma),\quad
&&\text{on}\ \partial \Omega\times (0,T),\\
&\mu_\Gamma=-\kappa \Delta_\Gamma \phi +  \partial_\mathbf{n}\phi + G'(\phi),\quad
&&\text{on}\ \partial \Omega\times (0,T),\\
&\phi|_{t=0}=\phi_0(x), \quad
&&\text{in}\ \Omega.
\end{aligned}
\right.
\label{CHbb}
\end{equation}
When the mobilities are positive constants, for instance, $M_{\mathrm{b}}=M_{\mathrm{s}}=1$, it was shown in \cite{GK20} that problem \eqref{CHbb} can be regarded as a $H^{-1}$ type gradient flow equation of the total free energy $E^{\mathrm{total}}(\phi)$ both in the bulk
and on the boundary, with respect to a suitable inner product on a certain function space (see \cite[Section 3]{GK20} for details). \medskip

\textbf{Well-posedness and long-time behavior}. Inspired by \cite{MZ05}, it will be convenient to view the trace of the order parameter $\phi$ as an unknown function on the boundary. After introducing the new variable
$$
\psi:=\phi|_{\partial \Omega},
$$
 the initial boundary value problem of the Cahn--Hilliard system \eqref{CHbb} can be written in the following form (taking $M_{\text{b}}=M_{\text{s}}=1$ for simplicity):
\begin{equation}
\left\{
\begin{aligned}
&\partial_t \phi=\Delta \mu,\quad
&&\text{in}\  \Omega\times (0,T),\\
&\mu=-\Delta \phi +F'(\phi),\quad
&&\text{in}\  \Omega\times (0,T),\\
&\partial_\mathbf{n} \mu=0,\quad
&&\text{on}\ \partial\Omega \times (0,T),\\
&\phi=\psi,\quad
&&\text{on}\ \partial\Omega \times (0,T),\\
&\partial_t \psi =\Delta_\Gamma \mu_\Gamma,\quad
&&\text{on}\ \partial\Omega \times (0,T),\\
&\mu_\Gamma= -\kappa \Delta_\Gamma \psi +\partial_\mathbf{n}\phi + G'(\psi),\quad
&&\text{on}\ \partial\Omega \times (0,T),\\
&\phi|_{t=0}=\phi_0(x), \quad
&&\text{in}\ \Omega,\\
&\psi|_{t=0}=\psi_0(x):=\phi_0(x)|_{\partial\Omega},\quad
&&\text{on}\ \partial\Omega .
\end{aligned}
\right.
\label{CH}
\end{equation}
Assume that $\Omega\subset \mathbb{R}^d$, $d=2,3$, is a smooth bounded domain.
Well-posedness of problem \eqref{CH} was first analyzed in \cite{LW} when $F$ and $G$ are
suitable regular potentials (including the typical choice like \eqref{regpo} for $F$ and the physically relevant surface potential in \cite{Q1} for $G$). The one dimensional case would be easier since the boundary dynamics is trivial as the terms involving $\Delta_\Gamma$ simply vanish.
The proof therein was inspired by the argument in \cite{MZ05}.
First, for the problem with surface diffusion (i.e., $\kappa > 0$), the authors introduced regularization
by adding viscous terms in both of the bulk and boundary chemical potentials.
This leads to a viscous Cahn--Hilliard equation subject to a dynamic boundary condition of viscous Cahn--Hilliard type,
which can be solved by using the contraction mapping principle (see \cite[Proposition 4.1]{LW}). After deriving uniform global-in-time a priori
estimates that are independent of the approximating parameter as well, they obtained the existence of global weak (and strong) solutions by
finding a convergent subsequence of the approximate solutions after passing to the limit, see \cite[Theorem 3.1]{LW}.

For the case without surface diffusion, i.e., $\kappa = 0$, the existence result can be achieved by deriving uniform estimates independent of the parameter $\kappa$ and then taking the limit as $\kappa \to 0^+$ (i.e., the vanishing surface diffusion limit), see \cite[Theorem 3.2]{LW}. In both cases, uniqueness of solutions can be proved by a standard energy method. It is worth mentioning that the solution will lose some spatial regularity due to the absence of the surface diffusion (cf. \cite{CF20b} for a similar situation for problem \eqref{claCHdb1}). In order to guarantee that the boundary equation $\mu_\Gamma= -\partial_\mathbf{n}\phi + G'(\psi)$ is satisfied for weak solutions  in the usual sense, that is, a.e. on $\partial\Omega \times (0,T)$, an additional geometric assumption was imposed in \cite[Theorem 3.2]{LW}. This restriction was later removed in \cite{GK20}, where the authors introduced a weaker notion of the ``weak solution" such that the chemical potentials satisfy the following weak form (see \cite[Definition 4.1]{GK20})
$$
\int_0^T\!\int_\Omega \mu \eta \,\d x\d t+\int_0^T\!\int_{\partial\Omega} \mu_\Gamma \eta\,\d S\d t =\int_0^T\!\int_\Omega  \nabla\phi\cdot \nabla \eta+F'(\phi)\eta \,\d x\d t +\int_0^T\!\int_{\partial\Omega}G'(\psi)\eta\,\d S\d t,
$$
for any test function $\eta\in L^2(0,T; H^1(\Omega))\cap L^\infty(\Omega\times (0,T))$ with $\eta|_{\partial\Omega}\in L^\infty(\partial\Omega\times(0,T))$. Namely, the troublesome boundary term $\partial_\mathbf{n} \phi$ does not appear in an explicit way. Then they proved existence and uniqueness of global weak solutions by an implicit time discretization based on the gradient flow structure of problem \eqref{CH}, see \cite[Theorem 4.3]{GK20}. Unfortunately, their argument does not work for singular potentials like the logarithmic potential \eqref{logpo} or the obstacle potential \eqref{obspo}.

Well-posedness of problem \eqref{CH} with singular potentials was first established in \cite{CFW20}. For the case $\kappa>0$, the authors treated the initial boundary value problem in a wide class of nonlinearities (covering \eqref{logpo}, \eqref{regpo} and \eqref{obspo}) with the compatibility condition that the boundary potential plays a dominating role. They proved existence and uniqueness of global weak solutions \cite[Theorem 2.3, Theorem 2.4]{CFW20} as well as existence of a unique global strong solution \cite[Theorem 4.1]{CFW20}. The proofs rely on several approximations of the original problem such as the introduction of viscous regularizations in the chemical potentials, the Yosida approximation for maximal monotone graphs and an implicit time discretization scheme for the bulk-boundary coupled system. They first proved the existence of a discrete solution by taking advantage of the general maximal monotone theory. After deriving a number of uniform estimates, they were able to conclude the existence results by performing limiting procedures with respect to the time step first and then the approximating parameters.

Concerning the long-time behavior of problem \eqref{CH}, the authors made a preliminary investigation for the case $\kappa>0$ in \cite{MW20} within the framework of infinite dimensional dynamical systems. For the system with regular potentials, they proved existence of exponential attractors, which also yields the existence of a global attractor with finite fractal dimension \cite[Theorem 2.1]{MW20}, while for the system with singular potentials, they showed the existence of a global attractor in a suitable complete metric space \cite[Theorem 2.2]{MW20}. The result is less satisfactory for the latter, because of the possible singularity of the bulk / boundary potentials at the pure phases $\pm 1$ and its interplay with the dynamic boundary condition (cf. \cite{MZ10}). To overcome this difficulty, one crucial step is to establish the strict separation property of the solution (cf. \cite{FW21} for problem \eqref{claCHdb1W}). Besides, the case without boundary diffusion ($\kappa=0$) turns out to be more difficult, since it is not clear so far whether the solutions to problem \eqref{CH} can define a $C_0$-semigroup in certain phase spaces.

The long-time behavior of a single solution was first analyzed in \cite{LW}. With the additional assumption that the regular potentials $F$ and $G$ are real analytic functions on $\mathbb{R}$, the authors proved that every global bounded weak/strong solution to problem \eqref{CH} will converge to a single equilibrium as $t\to+\infty$ and provided an estimate on the convergence rate, see \cite[Theorem 3.3]{LW}. The conclusion was achieved by employing an extended \L ojasiewicz--Simon type inequality \cite[Lemma 6.3]{LW}. We note that the above convergence results are valid for both cases with or without boundary diffusion, i.e., for $\kappa\geq 0$. Moreover, in presence of the boundary diffusion, by applying the \L ojasiewicz--Simon approach in a different way, the authors were able to give a further characterization on the Lyapunov stability of steady states (e.g., local energy minimizers) that may be non-isolated, see \cite[Theorem 3.4]{LW}. Whether the corresponding convergence or stability results hold for solutions to the system \eqref{CH} with singular potentials remains an open question, since it again relies on the regularity of solutions, in particular, the strict separation property from the pure states $\pm 1$.

At last, we remark that all the available results for problem \eqref{CH} were obtained for the case with constant mobilities. It will be interesting to study the case with non-constant mobilities that are concentration dependent or even degenerate.

\subsection{Some extended models}
In the last part, we mention some extensions of the Cahn--Hilliard equation with dynamic boundary conditions that have been studied in the recent literature.

We note that in the formulations \eqref{claCHdb1v2}, \eqref{claCHdb1Wg} and \eqref{CH}, some ``strong" relations were imposed on the phase-field function or the chemical potential via the trace operator, namely, in terms of some nonhomogeneous Dirichlet boundary conditions. To provide a more general description of the interactions between the materials in the bulk and the materials on the boundary, extended models with certain relaxed coupling relations between the bulk and boundary variables (e.g., via the Robin type boundary conditions) were introduced and analyzed in \cite{KL20,KLLM21,GKY22}. For an extended model with nonlocal effects, we refer to \cite{KS21}.

In \cite{KL20}, the authors considered the following bulk-boundary coupling system that can be regarded as an extension of \eqref{CH}:
\begin{equation}
\left\{
\begin{aligned}
&\partial_t \phi=\Delta \mu,\quad
&&\text{in}\  \Omega\times (0,T),\\
&\mu=-\Delta \phi +F'(\phi),\quad
&&\text{in}\  \Omega\times (0,T),\\
&\partial_\mathbf{n} \mu=0,\quad
&&\text{on}\ \partial\Omega \times (0,T),\\
&\phi=\alpha \psi + \beta,\quad
&&\text{on}\ \partial\Omega \times (0,T),\\
&\partial_t \psi =\Delta_\Gamma \mu_\Gamma,\quad
&&\text{on}\ \partial\Omega \times (0,T),\\
&\mu_\Gamma= -\kappa \Delta_\Gamma \psi + \alpha \partial_\mathbf{n}\phi + G'(\psi),\quad
&&\text{on}\ \partial\Omega \times (0,T),\\
&\phi|_{t=0}=\phi_0(x), \quad
&&\text{in}\ \Omega,\\
&\psi|_{t=0}=\psi_0(x):=\alpha^{-1}(\phi_0(x)|_{\partial\Omega}-\beta),\quad
&&\text{on}\ \partial\Omega,
\end{aligned}
\right.
\label{CHex}
\end{equation}
for some constants $\alpha\neq 0$, $\beta\in \mathbb{R}$. The affine linear transmission condition between the bulk and boundary phase-field variables (i.e., $\phi=\alpha \psi + \beta$ on $\partial\Omega$) was introduced in order to account for some non-trivial boundary interactions. This condition can be further relaxed by imposing a Robin type boundary condition:
\begin{equation}
\left\{
\begin{aligned}
&\partial_t \phi=\Delta \mu,\quad
&&\text{in}\  \Omega\times (0,T),\\
&\mu=-\Delta \phi +F'(\phi),\quad
&&\text{in}\  \Omega\times (0,T),\\
&\partial_\mathbf{n} \mu=0,\quad
&&\text{on}\ \partial\Omega \times (0,T),\\
&K\partial_\mathbf{n}\phi=H(\psi)-\phi,\quad
&&\text{on}\ \partial\Omega \times (0,T),\\
&\partial_t \psi =\Delta_\Gamma \mu_\Gamma,\quad
&&\text{on}\ \partial\Omega \times (0,T),\\
&\mu_\Gamma= -\kappa \Delta_\Gamma \psi + H'(\psi)\partial_\mathbf{n}\phi + G'(\psi),\quad
&&\text{on}\ \partial\Omega \times (0,T),\\
&\phi|_{t=0}=\phi_0(x), \quad
&&\text{in}\ \Omega,\\
&\psi|_{t=0}=\psi_0(x),\quad
&&\text{on}\ \partial\Omega,
\end{aligned}
\right.
\label{CHex2}
\end{equation}
for some $K>0$ and a function $H\in C^2(\mathbb{R})$. See also \cite{CFL19,LW20} for the analysis on some similar extended models of the Allen--Cahn equation with dynamic boundary conditions. Under suitable assumptions, the authors proved existence and uniqueness of global weak solutions to problems \eqref{CHex} and \eqref{CHex2} (see \cite[Theorem 2.1, Theorem 2.2]{KL20}). Besides, for the special case $H(s)=\alpha s+\beta$, they showed the weak convergence of solutions as the parameter $K\to 0^+$, and derived an error estimate between solutions of the two models \cite[Theorem 2.3]{KL20}.

Similar relaxations can be imposed on the chemical potentials, which account for possible reactions between the materials in $\Omega$ and on $\partial \Omega$. In \cite{KLLM21}, the authors studied the following problem
\begin{equation}
\left\{
\begin{aligned}
&\partial_t \phi= \Delta \mu,\quad &&\text{in}\ \Omega \times (0,T),\\
&\mu= - \Delta \phi + F'(\phi),\quad &&\text{in}\ \Omega \times (0,T),\\
&L\partial_\mathbf{n}\mu=\beta \mu_\Gamma-\mu && \text{on}\ \partial\Omega \times (0,T),\\
&\phi=\psi, && \text{on}\ \partial\Omega \times (0,T),\\
&\partial_t \psi +\beta \partial_\mathbf{n} \mu -  \Delta_\Gamma \mu_\Gamma =0, \quad && \text{on}\ \partial \Omega\times (0,T),\\
&\mu_\Gamma= - \kappa \Delta_\Gamma \psi +\partial_\mathbf{n}\phi +G'(\psi), \quad && \text{on}\ \partial \Omega\times (0,T),\\
&\phi|_{t=0}=\phi_0(x),\quad &&  \text{in}\ \Omega,\\
&\psi|_{t=0}=\psi_0(x)=\phi_0(x)|_{\partial\Omega},\quad && \text{on}\ \partial \Omega,
\end{aligned}
\right.
\label{claCHdb1Wga}
\end{equation}
for some constants $\beta\neq 0$ and $L>0$. The Robin type boundary condition indicates that the mass flux $\partial_\mathbf{n}\mu$ is driven by differences in the chemical potentials and the constant $L^{-1}$ can be interpreted as the reaction rate. Taking $\beta=1$, this boundary condition also allows people to establish a connection between models \eqref{claCHdb1Wg} and \eqref{CH}, via the formal limits $L\to 0^+$ and $L\to +\infty$, which correspond to the limit cases of an instantaneous reaction ($L\to 0^+$) and a vanishing reaction rate $L\to +\infty$, see \cite{KLLM21} for detailed discussions. Thus, problem \eqref{claCHdb1Wga} can be interpreted as an interpolation between \eqref{claCHdb1Wg} and \eqref{CH} where the $L$ corresponds to positive but finite kinetic rates. The authors proved the existence, uniqueness and regularity of weak solutions to problem \eqref{claCHdb1Wga} \cite[Theorem 3.1]{KLLM21} and investigated the asymptotic limits as $L\to 0^+$ and $L\to +\infty$, establishing also convergence rates for these limits, see \cite[Theorem 4.1, Theorem 4.2]{KLLM21}. Concerning long-time behavior of problem \eqref{claCHdb1Wga}, we refer to \cite{GKY22}, where the authors proved the existence of a global attractor and the convergence of global weak solutions to a single equilibrium as $t\to+\infty$ by means of a \L ojasiewicz--Simon type inequality \cite[Theorem 4.9, Theorem 4.13]{GKY22}. Besides, they proved that the global attractor of problem \eqref{claCHdb1Wg}
is stable with respect to perturbations of the kinetic rate and constructed a robust family exponential attractors for $L\in [0,1]$ (see \cite[Theorem 6.5]{GKY22}).

We note that the results obtained in \cite{KL20,KS21,KLLM21,GKY22} for those extended models are only valid for regular potentials including \eqref{regpo}. However, singular potentials like the logarithmic potential \eqref{logpo} or the obstacle potential \eqref{obspo} are not admissible. This will be an interesting field for future research.

\section*{Acknowledgements}
The work of the author was partially supported by NNSFC 12071084.


\end{document}